# The ternary operations of groupoids

## R. A. R. Monzo

**Abstract.** We investigate the ternary operations of groupoids and prove that there is a one-to-one correspondence between the collection of right modular groupoids with a left identity element $l$ and laterally commutative, $l$-bi-unital semiheaps. This result is applied to prove that the natural ternary operation induced by a right modular groupoid $S$ with left identity is isomorphic to the natural ternary operation of a groupoid $T$ if and only if $S$ and $T$ are isomorphic. Ternary operations of other classes of groupoids are also characterised, including those generalised heaps that appear as the natural or standard ternary operation of an inverse semigroup. The collections of such generalised heaps are proved to form algebraic varieties of type (3,1,1).

## 1. Introduction

In this paper we define, investigate and classify the "natural" and other ternary operations of particular classes of groupoids, the term "groupoid" being used to indicate a set with a binary operation. The ternary operations of semigroups were investigated by Hewitt and Zuckerman in 1968 [4]. Semiheaps are ternary operations satisfying the para-associative law and generalised heaps are idempotent semiheaps satisfying the two generalised heap axioms. Semiheaps and generalised heaps were first considered by Wagner [14]. More recently Lawson found a link between inverse semigroups and generalised heaps; namely, that "the Morita equivalence of inverse semigroups is determined by generalised heaps" [7]. The present work adds to that of these authors. Our particular interest is to investigate what an isomorphism between natural or "standard" ternary operations of groupoids implies about the structure of those groupoids, in particular whether it implies that they are isomorphic. We also explore what properties of ternary operations are preserved by ternary operation isomorphisms (Theorem 24).

In Section 2 we give definitions, establish some notation and extend an isomorphism theorem of Hewitt and Zuckerman. They proved that if $(S,\bullet)$ and $(S,\circ)$ are two semigroups on the set $S$, $(S,\bullet)$ is globally idempotent (that is, $S = S \bullet S$) and weakly reductive and the two (natural) ternary operations are equal, then $(S,\bullet)$ and $(S,\circ)$ are isomorphic [4, Theorem 22]. Inspired and guided by Hewitt and Zuckerman's proof of this result, we prove that if $S$ and $T$ are semigroups and either $S$ or $T$ is globally idempotent and weakly reductive, then isomorphism of the two natural ternary operations implies that $S$ and $T$ are isomorphic semigroups (Corollary 13). As a consequence, if the natural ternary operations of two semigroups $S$ and $T$ are isomorphic and either $S$ or $T$ has a right, left or two-sided identity then $S$ and $T$ are isomorphic (Corollary 14).

Also in Section 2, we prove that if a semigroup $S$ has an identity element and its natural ternary operation is isomorphic to that of a groupoid $T$, then if $T$ is medial or left lateral or if $x^2 = 1$ implies $x = 1$ in $S$, then $S$ and $T$ are isomorphic (Theorems 7 and 8). The same result holds if $S$ is a groupoid with identity and its natural ternary operation is a semiheap (Lemma 4). In addition, we prove that if $S$ is an idempotent groupoid with a right identity and $T$ is a distributive groupoid then $S$ and $T$ have isomorphic natural ternary operations if and only if they are isomorphic groupoids (Corollary 17.)

We prove that the natural ternary operation of a right modular groupoid is *not* necessarily a semiheap (Proposition 2, Section 2), but that if it either has a left identity or is left lateral or weakly associative, then the natural ternary operations of such right modular groupoids *are* laterally commutative semiheaps (Section 3, Theorem 18). In Section 3 we also establish a one-to-one correspondence between the collection of right modular groupoids with a left identity element and laterally commutative semiheaps with a bi-unitary element (Theorem 20). This extends Proposition 2 [3] to the following table:

| *Bi-unital Semiheap Type* | *Involuted monoid type* | *Groupoid type* |
|---|---|---|
| arbitrary | arbitrary | |
| idempotent semiheap | I- monoid | |
| generalised heap | inverse monoid | |
| heap | group | |
| laterally commutative | | right modular with a left identity |

Right modular groupoids with left identity are not listed in the second column in the table above because they are not involuted monoids. They do however have a unary operation that is an involutive automorphism and this operation is investigated in Section 6.

In Section 4 we prove that two right modular groupoids, each with a left identity element, have isomorphic natural ternary operations if and only if they are isomorphic (Theorem 22). The correspondence between laterally commutative, $l$-bi-unital semiheaps and right modular groupoids with a left identity $l$ is applied to prove that two laterally commutative, bi-unital semiheaps are isomorphic if and only if they naturally induce isomorphic right modular groupoids with left identities (Corollary 23).

In Section 5 we investigate "alpha-determined" groupoids, groupoids that are determined by an involutive, idempotent-fixed automorphism, alpha, on a semilattice of groups. The characterisation of these groupoids extends the result that a groupoid is a strongly regular (or completely inverse) $AG^{**}$ groupoid if and only if it is determined by an involutive automorphism on a semilattice of abelian groups [9]. We define a "standard" ternary operation on alpha-determined groupoids and prove that such ternary operations are generalised heaps. The standard generalised heaps of two such groupoids are proved to be isomorphic if and only if their underlying semilattices of groups are isomorphic.

In Section 6 we prove that in a weakly associative, right modular groupoid the natural ternary operation induced by a permutation on the set $\{1, 2, 3\}$ is a semiheap that also satisfies the generalised heap axioms. We prove that this is <u>not</u> the case in an arbitrary right modular groupoid with a left identity element. We also prove that in a right modular groupoid $S$ with a left identity element $l$, the unary operation $x \mapsto xl = x'$ satisfies three identities and that if a unary operation on $S$ satisfies those three identities, then it is the unary operation $x \mapsto xl$. Furthermore, we prove that if a groupoid with left identity element $l$ has a unary operation $'$ that satisfies these three identities, then $x' = xl$ and $[S]'$ is a semiheap if and only if $S$ is a right modular groupoid (Theorem 41). We prove also that the natural ternary operation induced by a right modular groupoid $S$ with left identity is isomorphic to the natural ternary operation of a groupoid $T$ if and only if $S$ and $T$ are isomorphic (Theorem 46). In this sense the structure of a right modular groupoid with left identity is completely determined by its natural ternary operation. Section 7 characterises types of units in ternary products of right modular groupoids with left identity. Section 8 characterises various ternary operations in certain classes of groupoids. Section 9 gives a new characterisation of heaps and explores the intimate connection between groups, Ward quasigroups, heaps and the natural ternary operations of Ward quasigroups. The last section, Section 10, characterises those generalised heaps that appear as the natural or standard ternary product of an inverse semigroup. A table summarising the results characterising ternary operations of various classes of groupoids is also given.

## 2. The natural ternary products of groupoids

Throughout this paper $xy \cdot z$ will denote the product $(xy)z$ and $x \cdot yz$ will denote the product $x(yz)$.

**Definition 1.** A groupoid is called **right modular** if $xy \cdot z = zy \cdot x$, **left lateral** if $x \cdot yz = y \cdot xz$ and **weakly associative** if $xy \cdot z = y \cdot xz$. Right modular, weakly associative groupoids are known as **$AG^{*}$ groupoids**. Right modular, left lateral groupoids are known as **$AG^{**}$ groupoids**.

A commutative, right modular groupoid is a commutative semigroup. A right modular groupoid with a right identity is commutative. Therefore, studies have focused on right modular groupoids with a left identity. Such an identity is unique in a right modular groupoid. Right modular groupoids are **medial** groupoids. That is, they satisfy the identity $xy \cdot zw = xz \cdot yw$. We will call a groupoid **reversible** if it satisfies the identity $xy \cdot zw = wz \cdot yx$. It follows easily that a right modular groupoid with a left identity is left lateral; that is, a right modular groupoid with a left identity is an $AG^{**}$ groupoid. A right modular groupoid with a left identity is also **paramedial**; that is, it satisfies the identity $xy \cdot zw = wy \cdot zx$. Consequently, a right modular groupoid with a left identity is reversible. Proposition 1 below, as far as the author is aware, does not appear in the literature.

**Proposition 1.** *If $S$ is a groupoid with left identity element $l$ then $S$ is a right modular groupoid if and only if it is paramedial or reversible.*

Proof. ($\Rightarrow$) As stated above, a right modular groupoid with left identity is paramedial.
($\Leftarrow$) If $S$ is a groupoid with left identity element $l$ and $S$ is paramedial or reversible then for all $a,b,c \in S$,
$ab \bullet c = ab \bullet lc \in \{cb \bullet la, cl \bullet ba\} = \{cb \bullet a\}$. ∎

**Definition 2.** Let $S$ be a set with a ternary operation [ ]. Then $\{S, [\ ]\}$, or simply $[S]$, is a **semiheap** if it satisfies the identities $[[xyz]\,uv] = [x\,[uzy]\,v] = [xy\,[zuv]]$. This is called the **para-associative law**. A ternary operation is **laterally commutative** if it satisfies the identity $[xyz] = [zyx]$. It is **left commutative** if it satisfies the identity $[xyz] = [yxz]$ and **right commutative** if it satisfies the identity $[xyz] = [xzy]$. It is **associative** if $[[xyz]uv] = [x[yzu]v] = [xy[zuv]]$.

We call $l \in S$ **left unital** if $[llx] = x$ for all $x \in S$ and **right unital** $[xll] = x$ for all $x \in S$. We call $l$ **bi-unitary in** $[S]$ if is is left and right unital. It is a **lateral unit** if $[lxl] = x$, for all $x \in S$ and an **outer lateral unit** if $[l\,[lxl]\,l] = x$, for all $x \in S$. A lateral unit is called a **unit** if it is bi-unitary. A ternary operation is **idempotent** if it satisfies the identity $[xxx] = x$, **l-bi-unital** if it has a bi-unitary element $l$, **l-unital** if it has a unit element $l$, and **left l-consistent** if it satisfies the identity $[xyz] = [l\,[lxy]\,z]$. An element $l$ is called a **central commutant** if $S$ satisfies the identity $[xly] = [ylx]$. The ternary operation $[S]$ is called **outer lateral** if every element of $S$ is outer lateral.

A semiheap is a **heap** if every element is biunitary. A semiheap is a **near heap** if it is idempotent and satisfies the identity $[xxy] = [yxx]$. A semiheap is a **generalised heap** if it is idempotent and satisfies the identities $[xx[yyz]] = [yy[xxz]]$ and $[[xyy]zz] = [[xzz]yy]$, called the **generalised heap axioms**. Two ternary operations $[S]$ and $\{T\}$ are **isomorphic** if there is a bijection $\delta : S \to T$ such that $\delta [xyz] = \{(\delta x)(\delta y)(\delta z)\}$ for all $x,y,z \in S$. In this case we write $[S] \cong \{T\}$ or $\delta : [S] \cong \{T\}$.

**Definition 3.** The **natural ternary operation** (or the **natural ternary product**) [ ] **induced by a groupoid** $S$ is defined as $[abc] = ab \bullet c$ for all $a,b,c \in S$. If $*$ is a unary operation on $S$ then the **ternary $*$ operation** $[\ ]*$ is defined as $[abc]* = ab^* \bullet c$, for all $a,b,c \in S$. Then $*$ is called **involutive** if $*^2 = 1_S$, where $1_S$ is the identity mapping on $S$, and **idempotent** if $*^2 = *$. Note that $[\ ] = [\ ]1_S$. A groupoid is called $*$**-unary** if it satisfies the identities, $(x^*)^* = x$, $(xy)^* = x^* y^*$ and $x^* y = y^* x$. In such a case, $*$ is an involutive automorphism and $(xy)^* = yx$. This implies that $i = i^*$ for any left or right identity in a $*$-unary groupoid and therefore left or right identities are unique. Note that a $*$-unary groupoid is reversible, since for all $x,y,z,w \in S$, $xy \bullet zw = (zw \bullet xy)^* = (zw)^* \bullet (xy)^* = wz \bullet yx$. $[S]$ is called $*$ **congruent** if it satisfies the identity $[xyz]^* = [x^* y^* z^*]$. A groupoid $S$ is called **outer lateral** if $[S]$ is outer lateral. It is obvious that any heap is outer lateral.

**Proposition 2.** The natural ternary operation on a right modular groupoid is not necessarily a semiheap.

Proof. The groupoid $S$ with the following multiplication table is right modular (cf. [8], page 88).

| S  | a  | b  | ab | ba |
|----|----|----|----|----|
| a  | a  | ab | ba | b  |
| b  | ba | b  | a  | ab |
| ab | b  | ba | ab | a  |
| ba | ab | a  | b  | ba |

Then $[[a(ba)b](ab)a] = [(a \bullet ba)b \bullet (ab)]\,a = (b \bullet ab)a = a \neq [a(ba)[b(ab)a]] = ba$. ∎

In the following section we will prove that if $S$ is a right modular, left lateral groupoid then $[S]$ is a semiheap.

**Proposition 3.** If $S \cong T$ then $[S] \cong [T]$.

Proof. Suppose that $\delta : S \cong T$. Let $a,b,c \in S$. Then $\delta [abc] = \delta (ab \bullet c) = (\delta a \delta b \bullet \delta c) = [(\delta a)(\delta b)(\delta c)]$. ∎

**Lemma 4.** *If* $\delta : [S] \cong [T]$ *and* $1 \in S$ *is the identity of S then*: (4.1) *for all* $a,b \in S$, $\delta(ab) = (\delta a \bullet \delta b)\delta 1$, (4.2) $(\delta 1)^2$ *is a left identity in T*, (4.3) *for all* $a \in S$, $\delta a = (\delta 1 \bullet \delta a)\delta 1 = (\delta a \bullet \delta 1)\delta 1$, (4.4) *for all* $a \in S$, $\delta a \bullet \delta 1 = \delta 1 \bullet \delta a$ *and* (4.5) $\delta : S \cong T$ *if and only if* $(\delta 1)^2 = \delta 1$.

*Furthermore, if [S] is a semiheap then S is commutative and if T is left lateral, medial or* $x^2 = 1$ *implies* $x = 1$ *in S or T then* $S \cong T$.

Proof. For all $a,b \in S$, $\delta(ab) = \delta(ab \bullet 1) = \delta[ab1] = [(\delta a)(\delta b)(\delta 1)] = (\delta a \bullet \delta b)\delta 1$, so (4.1) is valid. Also, $\delta a = \delta(11 \bullet a) = \delta[11a] = [(\delta 1)(\delta 1)(\delta a)] = (\delta 1)^2 \delta a$, so (4.2) is valid.
Then, $\delta a = \delta(1a \bullet 1) = \delta[1a1] = [(\delta 1)(\delta a)(\delta 1)] = (\delta 1 \bullet \delta a)\delta 1 = \delta(a1 \bullet 1) = = [(\delta a)(\delta 1)(\delta 1)] = = (\delta a \bullet \delta 1)\delta 1$, and so (4.3) is valid. Since $\delta a = (\delta 1 \bullet \delta a)\delta 1$, $\delta a \bullet \delta 1 = \{(\delta 1 \bullet \delta a) \delta 1\}\delta 1$. Now since $\delta$ is a surjection, there exists $x \in S$ such that $\delta x = \delta 1 \bullet \delta a$. Applying (4.3), $\delta a \bullet \delta 1 = \{(\delta 1 \bullet \delta a) \delta 1\}\delta 1 = = \{\delta x \bullet \delta 1\}\delta 1 = \delta x = \delta 1 \bullet \delta a$, so (4.4) is valid. If $\delta : S \cong T$ then $\delta 1$ is the identity of T. Since $(\delta 1)^2$ is a left identity of T, $\delta 1 = (\delta 1)^2$. Conversely, if $(\delta 1)^2 = \delta 1$ then by (2) and (4), $\delta a \bullet \delta 1 = \delta 1 \bullet \delta a = = (\delta 1)^2 \bullet \delta a = \delta a$. Then, by (4.1), $\delta(ab) = (\delta a \bullet \delta b)\delta 1 = (\delta a \bullet \delta b)$ and so $\delta : S \cong T$. So (4.5) is valid.

Suppose that [S] is a semiheap. Let $\delta u = (\delta 1)^2$. Then $\delta(u^2) = \delta(u^2 1) = \delta[uu1] = (\delta u)^2 \delta 1 = \delta 1$. So if $x^2 = 1$ implies $x = 1$ in S then $u^2 = 1$ imlies $u = 1$. Thus $\delta 1 = (\delta 1)^2$ and $S \cong T$ by 4.5. Now $\delta(ab) = = \delta[[1ab] 11] = \delta[1[1ba]1] = \delta(ba)$ and so $ab = ba$ and S is commutative. Then, by 4.2, $\delta a \bullet \delta b = = \delta u \delta a \bullet \delta b = \delta[uab] = \delta[uba] = \delta u \delta b \bullet \delta a = \delta b \bullet \delta a$. So T is commutative. By 4.2, $\delta u$ is the identity of T and, using (4.1), (4.2) and (4.3) and the commutivity of T it is straightforward to prove that if T is left lateral or medial then $\Phi : S \cong T$, where $\Phi x = (\delta x)(\delta 1)$. ∎

**Lemma 5.** *If* $\delta : [S] \cong [T]$ *and* $1 \in S$ *is the identity element of S then if T is medial, left lateral or a semigroup then* $(\delta 1)^2$ *is the identity in T and for all* $a,b,c \in S$, $\delta(a \bullet bc) = \delta a \bullet \delta b \delta c$.

Proof. Let $a \in S$. Using (4.3), $\delta a = (\delta a \bullet \delta 1)\delta 1$. If T is a semigroup then $\delta a = (\delta a \bullet \delta 1)\delta 1 = = (\delta a)(\delta 1)^2$ and, by (4.2), $(\delta 1)^2$ is an identity element. If T is medial then, using (4.2) and (4.4), $\delta a = = (\delta a \bullet \delta 1)\delta 1 = \delta 1(\delta a \bullet \delta 1) = \{(\delta 1)^2 \bullet \delta 1\}(\delta a \bullet \delta 1) = \{(\delta 1)^2 \bullet \delta a\}(\delta 1)^2 = \delta a (\delta 1)^2$ and so $(\delta 1)^2$ is an identity element. If T is left lateral then using (4.4), $\delta a = (\delta a \bullet \delta 1)\delta 1 = \delta 1(\delta a \bullet \delta 1) = (\delta a)(\delta 1)^2$ and so $(\delta 1)^2$ is an identity element.

Now, using (4.1), (4.3) and (4.4), for all $a,b,c \in S$, $\delta(a \bullet bc) = \{\delta a \bullet \delta(bc)\}\delta 1 = \delta 1\{\delta a \bullet \delta(bc)\}$ and if T is medial or satisfies the identity $x \bullet yz = y \bullet xz$ then $\delta(a \bullet bc) = \delta a\{\delta 1 \bullet \delta(bc)\} = \delta a\{\delta(bc) \bullet \delta 1\} = = \delta a \bullet \{(\delta b \delta c) \bullet \delta 1\}\delta 1 = \delta a \bullet \delta b \delta c$.

If T is a semigroup then $\delta(a \bullet bc) = \{\delta a \bullet \delta(bc)\}\delta 1 = \delta a\{\delta(bc) \bullet \delta 1\} = \delta a \bullet \{(\delta b \delta c) \bullet \delta 1\}\delta 1 = = \delta a \bullet \delta b \delta c$. ∎

**Theorem 6.** *Suppose that* $\delta : [S] \cong [T]$ *and* $1 \in S$ *is the identity element of S . If T is right distributive then* $\delta : S \cong T$ *and S is an idempotent groupoid.*

Proof. Since T is right distributive, by (4.2) and (4.3), for all $a \in S$, $\delta 1 = (\delta 1)^2 \delta 1 = (\delta 1)^2 (\delta 1)^2 = (\delta 1)^2$. Then by (4.4), $\delta a = (\delta 1)^2 \bullet \delta a = \delta 1 \bullet \delta a = \delta a \bullet \delta 1$ and so $\delta 1 = (\delta 1)^2$ is an identity element. Using (4.1), for all $a,b,c \in S$, $\delta(a \bullet bc) = \{\delta a \bullet \delta(bc)\}\delta 1 = \delta a \bullet \delta(bc) = \delta a \bullet \{(\delta b \delta c) \bullet \delta 1\} = \delta a \bullet (\delta b \delta c)$. Putting $c = 1$ gives $\delta(ab) = \delta a \bullet (\delta b \delta 1) = \delta a \bullet \delta b$ and so $\delta : S \cong T$.

Note that $\delta a = (\delta 1)^2 \bullet \delta a = (\delta 1 \bullet \delta a)^2 = (\delta a)^2$. Hence, T (and therefore S) is idempotent. ∎

**Theorem 7.** *Suppose that* $\delta : [S] \cong [T]$, $1 \in S$ *is the identity element of S and S is a semigroup. If T is medial, left lateral or a semigroup then* $S \cong T$.

Proof. Since $S$ is a semigroup, $\delta(ab \bullet c) = \delta(a \bullet bc)$. Hence, using Lemma 5 and the hypothesis that $\delta: [S] \cong [T]$, $\delta[abc] = \delta(ab \bullet c) = \delta(a \bullet bc) = [(\delta a)(\delta b)(\delta c)] = \delta a \delta b \bullet \delta c = \delta a \bullet \delta b \delta c$ and so $T$ is a semigroup. Let $a,b \in S$. Define $\Phi a = \delta a \delta 1$. Then, using (4.1) and (4.3), $\Phi(ab) = \delta(ab)\delta 1 = (\delta a \delta b \bullet \delta 1)\delta 1 = \delta a \delta b = \delta a (\delta 1 \delta b \bullet \delta 1) = (\delta a \delta 1)(\delta b \delta 1) = \Phi a \Phi b$. If $\Phi a = \Phi b$ then $\delta a \delta 1 = \delta b \delta 1$ and so $\delta a = (\delta a \delta 1)\delta 1 = (\delta b \delta 1)\delta 1 = \delta b$, so $\Phi$ is injective. If $\delta x = \delta a \delta 1$ then $\Phi x = \delta x \delta 1 = (\delta a \delta 1)\delta 1 = \delta a$, so $\Phi$ is surjective. Therefore, $\Phi: S \cong T$. ∎

**Theorem 8.** *If $\delta: [S] \cong [T]$, $1 \in S$ is the identity element of $S$ and $x^2 = 1$ implies $x = 1$ in $S$, then $S \cong T$.*

Proof: Let $a,b \in S$ and let $\delta x = (\delta 1)^2$. Then, using Lemma 4, $\delta(x^2 a) = \delta(x^2 a \bullet 1) = \delta[x^2 a 1] = [(\delta x^2)(\delta a)(\delta 1)] = (\delta x^2 \delta a)\delta 1 = \{(\delta x)^2 \delta 1 \bullet \delta a\}\delta 1 = (\delta 1 \bullet \delta a)\delta 1 = \delta a$ and so $a = x^2 a$. Hence, $1 = x^2 1 = x^2$, which implies $x = 1$. Therefore, $\delta 1 = (\delta 1)^2$ and by (4.5), $S \cong T$. ∎

**Corollary 9.** *Let $S$ be a semigroup without an identity. Then $[S^1] \cong [T]$ if and only if $S^1 \cong T$, where $S^1$ is the semigroup formed by adjoining an identity element 1 to $S$.*

**Definition 4.** A groupoid $S$ is **weakly reductive** if for all $x,y \in S$, $xz = yz$ and $zx = zy$ for all $z \in S$ implies $x = y$. It is **globally idempotent** if $S = S^2$.

**Lemma 10.** *If $\delta: [S] \cong [T]$ and $S$ is globally idempotent then $T$ is globally idempotent.*

Proof. For $\delta x \in T$, since $x = ab = (cd)b$ for some $a,b,c,d \in S$, $\delta x = \delta(cd \bullet b) = \delta[cdb] = \delta c \delta d \bullet \delta b \in T^2$. ∎

**Lemma 11.** *If $\delta: [S] \cong [T]$, $S$ and $T$ are semigroups and $S$ is globally idempotent and weakly reductive then $T$ is weakly reductive.*

Proof. Consider any $\delta x, \delta y \in T$. Suppose that $\delta z \delta x = \delta z \delta y$ and $\delta x \delta z = \delta y \delta z$ for all $z \in S$. Now, for any $a,b \in S$, $\delta a \delta b = \delta c$ for some $c \in S$. Hence, $\delta c \delta x = \delta c \delta y$ and $\delta x \delta c = \delta y \delta c$. Therefore, $\delta(ab \bullet x) = \delta[abx] = [(\delta a)(\delta b)(\delta x)] = \delta a \delta b \bullet \delta x = \delta a \delta b \bullet \delta y = [(\delta a)(\delta b)(\delta y)] = \delta[aby] = \delta(ab \bullet y)$ and so $ab \bullet x = ab \bullet y$. Since $S = S^2$, for any $w \in S$, $wx = wy$.

Also, $\delta x (\delta a \delta b) = \delta y (\delta a \delta b)$. But since $S$ and $T$ are semigroups, $\delta(xab) = \delta[xab] = (\delta x \delta a)\delta b = \delta x (\delta a \delta b) = \delta y (\delta a \delta b) = (\delta y \delta a)\delta b = \delta[yab] = \delta(yab)$ and so $xab = yab$. Therefore, $xw = yw$, for any $w \in S$. Then since $S$ is weakly reductive, $x = y$. and so $\delta x = \delta y$ and $T$ is weakly reductive. ∎

**Theorem 12.** *If $\delta: [S] \cong [T]$, $S$ and $T$ are semigroups, $S$ is globally idempotent and weakly reductive then $S \cong T$.*

Proof. For each $x \in S$ we choose fixed elements $l_x$ and $r_x$ such that $x = l_x r_x$. We will show that the mapping $\Phi: S \to T$; $x \mapsto \delta l_x \delta r_x$ is an isomorphism. If $\Phi x = \Phi y$ then for all $z \in S$, $\delta(xz) = \delta(l_x r_x \bullet z) = \delta[l_x r_x z] = [\delta l_x \delta r_x \delta z] = \delta l_x \delta r_x \bullet \delta z = \delta l_y \delta r_y \bullet \delta z = \delta[l_y r_y z] = \delta(l_y r_y \bullet z) = \delta(yz)$. Hence, $xz = yz$. Then, $\delta(zx) = \delta(z \bullet l_x r_x) = \delta(z l_x \bullet r_x) = \delta[z l_x r_x] = [\delta z \delta l_x \delta r_x] = \delta z \delta l_x \bullet \delta r_x = \delta z \bullet \delta l_x \delta r_x = \delta z \bullet \delta l_y \delta r_y = \delta z \delta l_y \bullet \delta r_y = \delta[z l_y r_y] = \delta(z l_y \bullet r_y) = \delta(z \bullet l_y r_y) = \delta(zy)$. Hence, $zx = zy$. Since $S$ is weakly reductive, $x = y$ and $\Phi$ is injective.

Suppose that $\delta x \in T$. Then, by Lemma 10, $\delta x = \delta v \delta w$ for some $v,w \in S$. Then, for any $z \in S$, $\delta x \delta z = \delta v \delta w \bullet \delta z = [(\delta v)(\delta w)(\delta z)] = \delta[vwz] = \delta(vw \bullet z) = \delta(l_{vw} r_{vw} \bullet z) = \delta[l_{vw} r_{vw} z] = \delta l_{vw} \delta r_{vw} \bullet \delta z = \Phi(vw) \delta z$.

Then $\delta z \delta x = \delta z \bullet \delta v \delta w = \delta z \delta v \bullet \delta w = [(\delta z)(\delta v)(\delta w)] = \delta [zvw] = \delta (zv\bullet w) = \delta (z\bullet vw) =$
$= \delta (z\bullet l_{vw} r_{vw}) = \delta (zl_{vw} \bullet r_{vw}) = \delta [zl_{vw} r_{vw}] = \delta z \delta l_{vw} \bullet \delta r_{vw} = \delta z \bullet \delta l_{vw} \delta r_{vw} = \delta z \Phi(vw)$.
Now since by Lemma 11, $T$ is weakly reductive, $\delta x = \Phi(vw)$ and $\Phi$ is surjective.

We have proved that $\Phi$ is bijective and we need only show that it is a homomorphism. Let $x,y,z \in S$.
Then $\delta z \Phi(xy) = \delta z \bullet \delta l_{xy} \delta r_{xy} = \delta z \delta l_{xy} \bullet \delta r_{xy} = [(\delta z)(\delta l_{xy})(\delta r_{xy})] = \delta [z l_{xy} r_{xy}] =$
$= \delta (z l_{xy} \bullet r_{xy}) = \delta (z \bullet l_{xy} r_{xy}) = \delta(z\bullet xy)$. Also, $\Phi(xy) \delta z = \delta l_{xy} \delta r_{xy} \bullet \delta z =$
$= [(\delta l_{xy})(\delta r_{xy})(\delta z)] = \delta [l_{xy} r_{xy} z] = \delta (l_{xy} r_{xy} \bullet z) = \delta(xy\bullet z)$.

Then, $\Phi x \delta y = \delta l_x \delta r_x \bullet \delta y = [(\delta l_x)(\delta r_x)(\delta y)] = \delta [l_x r_x y] = \delta (l_x r_x \bullet y) = \delta(xy)$.

Also, $\delta y \Phi x = \delta y \bullet \delta l_x \delta r_x = \delta y \delta l_x \bullet \delta r_x = [(\delta y)(\delta l_x)(\delta r_x)] = \delta [y l_x r_x] = \delta (yl_x \bullet r_x)$
$= \delta (y\bullet l_x r_x) = \delta(yx)$.

So $\delta z(\Phi x \Phi y) = (\delta z \Phi x) \Phi y = \delta(zx) \Phi y = \delta (zx\bullet y) = \delta(z\bullet xy) = \delta z \Phi(xy)$.
Also, $(\Phi x \Phi y) \delta z = \Phi x (\Phi y \delta z) = \Phi x \delta(yz) = \delta(x\bullet yz) = \delta(xy\bullet z) = \Phi(xy) \delta z$. Then, since by Lemma 11 again, $T$ is weakly reductive, $\Phi(xy) = \Phi x \Phi y$ and $\Phi$ is a homomorphism. Hence, $\Phi : S \cong T$. ∎

**Corollary 13.** *If $\delta : [S] \cong [T]$, $S$ and $T$ are semigroups and either $S$ or $T$ is globally idempotent and weakly reductive then $S \cong T$.*

Proof. This follows because $\delta^{-1}[(\delta x)(\delta y)(\delta z)] = \delta^{-1} \delta [xyz] = [xyz] =$
$= [(\delta^{-1} \delta x)(\delta^{-1} \delta y)(\delta^{-1} \delta z)]$ and, since $\delta$ and $\delta^{-1}$ are both bijections, $\delta^{-1} : [T] \cong [S]$. ∎

**Corollary 14.** *If $\delta : [S] \cong [T]$, $S$ and $T$ are semigroups and either $S$ or $T$ has a right, left, or two-sided identity then $S \cong T$.*

Proof. If a groupoid has a right, left or two-sided identity then it is globally idempotent and weakly reductive and so Corollary 14 follows from Corollary 13. ∎

**Lemma 15.** *If $\delta : [S] \cong [T]$, $S$ has a right identity element $r \in S$ and $T$ is right distributive then $S$ is right distributive. If $T$ is also left distributive then $T$ is idempotent.*

Proof. If $a,b \in S$ then $\delta(ab) = \delta(ab\bullet r) = \delta[abr] = \delta a \delta b \bullet \delta r$ and so $\delta a = \delta a \delta r \bullet \delta r$. This implies that $\delta(ab) \delta r = \delta a \delta b$. Then, using the right distributivity of T, for any $a,b,c \in S$, $\delta(ab\bullet c) = \delta[abc] =$
$= \delta a \delta b \bullet \delta c = \delta a \delta c \bullet \delta b \delta c = \delta(ac) \delta r \bullet \delta(bc) \delta r = \delta(ac) \delta(bc) \bullet \delta r = \delta(ac\bullet bc)$ and so $ab\bullet c = ac \bullet bc$. That is, $S$ is right distributive.

Now we assume that $T$ is distributive. We have $\delta r = \delta r \delta r \bullet \delta r = (\delta r)^2 (\delta r)^2 = (\delta r)^2 \bullet (\delta r)^2 (\delta r)^2 =$
$= (\delta r)(\delta r \bullet \delta r \delta r) = (\delta r)\{(\delta r)^2 \bullet (\delta r)^2\} = (\delta r)^2$. Then, $\delta a = \delta a \delta r \bullet \delta r = \delta a (\delta r)^2 \bullet \delta r =$
$= (\delta a \delta r)^2 \bullet \delta r = (\delta a \delta r \bullet \delta r)^2 = (\delta a)^2$ and so $T$ is idempotent. ∎

**Theorem 16.** *If $\delta : [S] \cong [T]$, $S$ is idempotent and has a right identity element $r \in S$ and $T$ is distributive then $S \cong T$.*

Proof. As in the proof of Lemma 15, for all $a,b \in S$, $\delta a = \delta a^2 = (\delta a)^2 \delta r = \delta a \delta r$. Then, as in the proof of Lemma 15 again, $\delta(ab) = \delta a \delta b \bullet \delta r = \delta a \delta b$ and so $\delta : S \cong T$. ∎

**Corollary 17.** *If $S$ is idempotent and has a right identity element $r \in S$ and $T$ is distributive then $[S] \cong [T]$ if and only if $S \cong T$.*

Proof. This follows from Theorem 16 and Proposition 3. ∎

## 3. Laterally commutative semiheaps and right modular groupoids

**Note 1.** If $S$ is an $AG^*$ groupoid then it satisfies the **permutation identity** $x_1 x_2 \cdot x_3 x_4 = x_{\Pi(1)} x_{\Pi(2)} \cdot x_{\Pi(3)} x_{\Pi(4)}$, where $\Pi$ is any permutation of $\{1, 2, 3, 4\}$ [12]. An $AG^*$ groupoid also satisfies the identity $xy \cdot z = y \cdot xz = y \cdot zx$.

**Theorem 18.** *If $S$ is an $AG^*$ or an $AG^{**}$ groupoid then $[S]$ is a laterally commutative semiheap.*

Proof. Let $S$ be an $AG^{**}$ groupoid. For any $a, b, c \in S$, $[[abc]\,de] = \{(ab \cdot c)d\}e = (dc \cdot ab)\,e = \{a(dc \cdot b)\}e = [a\,[dcb]\,e] = \{a(bc \cdot d)\}e = \{e(bc \cdot d)\}a = (bc \cdot ed)a = (a \cdot ed)(bc) = (ab)(ed \cdot c) = (ab)(cd \cdot e) = [ab\,[cde]\,]$. So $[S]$ is a semiheap. Clearly, since $S$ is right modular, $[S]$ is laterally commutative.

Let $S$ be an $AG^*$ groupoid. Then $S$ satisfies the permutation identity $x_1 x_2 \cdot x_3 x_4 = x_{\Pi(1)} x_{\Pi(2)} \cdot x_{\Pi(3)} x_{\Pi(4)}$, where $\Pi$ is any permutation of $\{1,2,3,4\}$ [11]. So $[[abc]\,de] = \{(ab \cdot c)d\}e = (dc \cdot ab)e = (cd \cdot ab)e = (ab)(cd \cdot e) = [ab\,[cde]]$. Also, $[a\,[dcb]\,e] = a(dc \cdot b) \cdot e = (dc \cdot b)(ae) = (bc \cdot d)(ae) = d(bc \cdot ae) = d(ab \cdot ec) = (ab)(d \cdot ec) = (ab)(cd \cdot e) = [ab\,[cde]]$. So $[S]$ is a semiheap. Since $S$ is right modular, $[S]$ is laterally commutative. ∎

As mentioned above, a right modular groupoid $S$ with a left identity element is an $AG^{**}$ groupoid. It follows from Theorem 18 then, that $[S]$ is a laterally commutative semiheap. However, we shall use the following Lemma to prove that $[S]$ is a semiheap when $S$ is a right modular groupoid with a left identity.

**Lemma 19.** *If $S$ has a laterally commutative ternary product $[\ ]$ with a bi-unitary element $l$ then $[S]$ is a semiheap if and only if it satisfies the identity $[[abc]\,de] = [[cdl]\,[bal]\,e]$.*

Proof. ($\Rightarrow$) $[\,[cdl]\,[bal]\,e\,] = [\,[ll\,[cdl]]\,[bal]\,e\,] = [\,l\,[[bal]\,[cdl]\,l]\,e\,] = [\,e\,[[bal]\,[cdl]\,l]\,l\,] = [\,[el\,[cdl]]\,[bal]\,l\,] = [l\,[lab]\,[el\,[cdl]]\,] = [\,[lla]\,b\,[el\,[cdl]]\,] = [\,a\,b\,[el\,[cdl]]\,] = [\,a\,b\,[[cdl]\,le]\,] = [\,a\,b\,[cd\,[lle]\,]= [ab\,[cde]\,] = [[abc]\,de]$

($\Leftarrow$) We have **(1)**: $[[abc]\,de] = [[cdl]\,[bal]\,e]$. Therefore $[abc] = [[abc]\,ll]$ and it follows from (1) that $[abc] = [\,[cll]\,[bal]\,l\,] = [c\,[bal]\,l]$. So **(2):** $[abc] = [c\,[bal]\,l]$. From (1) it also follows that $a = [\,[all]\,ll\,] = [\,[lll]\,[lal]\,l\,] = [\,l\,[lal]\,l\,]$ and so we have **(3):** $a = [\,l\,[lal]\,l\,]$. Finally, from (2) we have **(4):** $[bal] = [l\,[abl]\,l]$. Then from (1),

$$\begin{aligned}
[[abc]\,de] &= [[cdl]\,[bal]\,e] = [e\,[[bal]\,[cdl]\,l]\,l] &&\text{(by (2))} \\
&= [e\,[[dcb]\,al]\,l\,] &&\text{(by (1))} \\
&= [l\,[[dcb]\,al]\,e\,] &&\text{(by lateral commutivity)} \\
&= [[lll]\,[[dcb]\,al]\,e\,] &&\text{($l$ is bi-lateral)} \\
&= [\,[a\,[dcb]\,l]\,le\,] &&\text{(by (1))} \\
&= [\,[a\,[dcb]\,l]\,[lll]e\,] &&\text{($l$ is bi-lateral)} \\
&= [\,[lla]\,[dcb]\,e\,] &&\text{(by (1))} \\
&= [a\,[dcb]\,e] &&\text{($l$ is bi-lateral)}
\end{aligned}$$

So we have **(5):** $[[abc]\,de] = [a\,[dcb]\,e]$

Also, $[a\,[dcb]\,e] = [\,[l\,[lal]\,l]\,[b\,[cdl]\,l]\,e\,]$    (by (2) and (3))
$\qquad = [\,[[cdl]\,bl]\,[lal]\,e]$    (by (1))
$\qquad = [\,[al\,[cdl]]\,be]$    (by (1))
$\qquad = [\,[[cdl]\,la]\,be]$    (lateral commutivity)
$\qquad = [\,[[cdl]\,[lll]a]\,be]$    ($l$ is bi-lateral)
$\qquad = [\,[[llc]\,da]\,be]$    (by (1))
$\qquad = [\,[cda]\,be]$    ($l$ is bi-lateral)
$\qquad = [\,eb\,[cda]]$    (lateral commutivity)

Therefore, $[a\,[dcb]\,e] = [e\,[dcb]\,a] = [ab\,[cde]\,]$. So we have proved that $[[abc]\,de] = [a\,[dcb]\,e] = [ab\,[cde]]$. Hence, $[S]$ is a semiheap. This completes the proof of Lemma 2. ∎

**Theorem 20.** *Let S be a right modular groupoid with a left identity element l and let [ ] be the natural ternary operation on S. Then [S] is a laterally commutative l-bi-unital semiheap. Conversely, if S is a laterally commutative, l-bi-unital semiheap then it induces a right modular groupoid with left identity element l if we define a product $*$ on S by $a*b = [lab]$.*

*Furthermore, if $\Phi$ maps the right modular groupoid S with left identity element l to the laterally commutative, l-bi-unital semiheap [S] and $\Gamma$ maps the laterally commutative semiheap S with bi-unitary element l to $\Gamma(S) = \{S, *, l\}$, as defined above, then $\Phi$ and $\Gamma$ are mutually inverse mappings.*

Proof. Let $(S, l)$ be a right modular groupoid with a left identity element $l$ and let $[S]$ be the natural ternary operation on $S$. Then, since $[abc] = (ab)c = (cb)a = [cba]$, $[S]$ is laterally commutative. Also, since $[lla] = (ll)a = la = a = (al)l = [all]$, $l$ is bi-unitary. Furthermore, $[[cdl][bal]e] = \{(cd\bullet l)(ba\bullet l)\}e =$
$= \{(cd\bullet ba)l\}e = \{(ba\bullet cd)\}e = (dc\bullet ab)e = \{(ab\bullet c)d\}e = [[abc]de]$. Therefore, by Lemma 19, $[S]\}$ is a semiheap. Hence, $[S]$ is a laterally commutative, $l$-bi-unital semiheap. We denote $[S]$ by $\Phi(S, l)$ or $[(S, l)]$.

If $\{S,[ ]\}$, or $[S, l]$, is a laterally commutative, $l$-bi-unital semiheap, define a product on $S$ by $a*b = [lab]$. Then $l*a = [lla] = a$ and so $l$ is a left identity element. Also, $(a*b)*c = [l(ab)c] = [l[lab]c] =$
$= [l[bal]c] = [[lla]bc]] = [abc] = [cba] = (c*b)*a$. Hence, this product produces a right modular groupoid structure on $S$, with left identity element $l$. We denote this groupoid as $\Gamma[S, l]$, $([S, l])$ or $\{S, *\}$.

Note that in $\Phi\Gamma[S, l]$, $[abc] = (a*b)*c = [l[lab]c] = [c[lab]l] = [[cba]ll] = [cba] = [abc]$ and so $\Phi\Gamma[S, l] = [S, l]$. Also, in $\Gamma\Phi(S, l)$, $a*b = [lab] = (la)b = ab$ and so $\Gamma\Phi(S, l) = (S, l)$. This completes the proof of Theorem 20. ∎

**Definition 5.** We call $([S, l])$ the ***natural right modular groupoid with left identity induced by the l-bi-unitary, laterally commutative semiheap*** $[S, l]$. We call $[(S, l)]$ ***the natural l-bi-unitary, laterally commutative semiheap induced by the right modular groupoid with left identity*** $(S, l)$.

## 4. Isomorphism theorems

Groupoids with isomorphic natural ternary operations may not be isomorphic. In fact, such groupoids exist of every order greater than or equal to 2.

**Example 1.** Let $S = \{x_\alpha : \alpha \in \Omega\}$ be a set of order $\Omega \geq 2$. Let $\{S, *\}$ be the left zero semigroup on $S$. That is, $x_\alpha * x_\beta = x_\alpha$ ($\alpha, \beta \in \Omega$). Fix $\alpha$ and $\beta$ in $\Omega$. Define $\{S, \bullet\}$ as follows: $x_\alpha \bullet x_\lambda = x_\beta$ ($\lambda \in \Omega$), $x_\beta \bullet x_\lambda = x_\alpha$ ($\lambda \in \Omega$) and $x_\gamma \bullet x_\lambda = x_\gamma$ ($\gamma, \lambda \in \Omega$ and $\gamma \notin \{\alpha, \beta\}$). Let $\theta = 1_S$. Then it is straightforward to show that $\theta: [\{S, *\}] \cong [\{S, \bullet\}]$, while clearly $\{S, *\}$ and $\{S, \bullet\}$ are not isomorphic.

We proceed to prove that right modular groupoids with a left identity are "determined" by their naturally induced ternary operations, in the sense that isomorphism between the ternary operations of two right modular groupoids with left identity implies isomorphism of the groupoids. Similarly, we prove that $l$-bi-unital, laterally commutative semiheaps are "determined" by the natural right modular groupoids (with left identity element $l$) that they induce. Later (Theorem 46) we will prove that, unlike in the Example 1 above, if $\{S, l\}$ is a right modular groupoid with left identity element $l$ then $[S] \cong [T]$ if and only if $S \cong T$.

**Lemma 21.** *Let $(S, l)$ and $(T, k)$ be right modular groupoids with left identities l and k respectively. Then $[(S, l)] \cong [(T, k)]$ if and only if there is a bijection $\delta: S \to T$ such that for all $a, b \in S$, $\delta(ab) = (\delta l \bullet \delta a)\delta b$.*

Proof. ($\Rightarrow$) Let $\delta: [(S, l)] \cong [(T, k)]$. Then for $a, b \in S$, $ab = ba\bullet l = [bal]$ and so $\delta(ab) = \delta[bal] = [(\delta b)(\delta a)(\delta l)] = (\delta b \bullet \delta a)\delta l = (\delta l \bullet \delta a)\delta b$.

($\Leftarrow$) Note that since for all $a, b \in S$, $\delta(ab) = (\delta l \bullet \delta a)\delta b$, $\delta b = \delta(lb) = (\delta l)^2 \delta b$ and so $(\delta l)^2 = k$. So $\delta[abc] = \delta(ab\bullet c) = (\delta l \bullet \delta(ab))\delta c = \delta l\{(\delta l \bullet \delta a)\delta b\}\bullet (\delta l)^2 \delta c = \delta l(\delta l)^2 \bullet \{(\delta c \bullet \delta b)(\delta l \bullet \delta a)\} =$

= { $\delta l (\delta c \bullet \delta b)$ }($\delta l \bullet \delta a$) = $(\delta l)^2$ {$(\delta c \bullet \delta b) \delta a$} = $(\delta a \bullet \delta b) \delta c$ = [$(\delta a)(\delta b)(\delta c)$] and therefore $\delta$ : [(S, l)] $\cong$ [(T, k)]. This completes the proof of Lemma 21.  ∎

**Theorem 22.** *Let (S, l) and (T, k) be right modular groupoids with left identities l and k respectively. Then (S, l) $\cong$ (T, k) if and only if* **[(S, l)] $\cong$ [(T, k)].**

Proof. ($\Rightarrow$) This follows from Proposition 3. ($\Leftarrow$) If **[(S, l)] $\cong$ [(T, k)]** then, by Lemma 21, there exists a bijection $\delta : S \to T$ *such that for all a,b $\in$ S*, $\delta (ab) = (\delta l \bullet \delta a) \delta b$ and so $(\delta l)^2 = k$.

Define $\Phi : S \to T$ as $\Phi a = \delta l \bullet \delta a$ ($a \in S$). We will prove that $\Phi : (S, l) \cong (T, k)$. Firstly, $\Phi(ab) = \delta l \bullet \delta (ab) =$
= $(\delta l)^2 \delta l \bullet \{(\delta l \bullet \delta a) \delta b\} = (\delta l)^2 \delta l \bullet \{(\delta b \bullet \delta a) \delta l\} = (\delta b \bullet \delta a) (\delta l)^2 = (\delta l)^2 \delta a \bullet \delta b = \delta a \bullet \delta b =$
= $(\delta l)^2 (\delta a \bullet \delta b) = (\delta l \bullet \delta a)(\delta l \bullet \delta b) = \Phi a \Phi b$. Therefore, $\Phi$ is a homomorphism.

Now $\delta l \bullet \delta a = (\delta a \bullet \delta l)(\delta l)^2$. So $(\delta l \bullet \delta a)(\delta l)^2 = (\delta a \bullet \delta l)$ and $\{(\delta l \bullet \delta a)(\delta l)^2\} \delta l = \delta a$. Hence, $\Phi(a) = \Phi(b)$ implies $\delta l \bullet \delta a = \delta l \bullet \delta b$ and then $\delta a = \{(\delta l \bullet \delta a)(\delta l)^2\} \delta l = \{(\delta l \bullet \delta b)(\delta l)^2\} \delta l = \delta b$. Since $\delta$ is a bijection, $a = b$ and $\Phi$ is an injection.

Let $a \in S$. Then $\Phi(al) = \delta l \bullet \delta (al) = (\delta l)^2 \delta l \bullet \{(\delta l \bullet \delta a) \delta l\} = (\delta l \bullet \delta a)(\delta l)^2 = (\delta a \bullet \delta l)$.
Now, for any $a \in S$ we define $a^*$ as follows: $\delta a^* = (\delta a \bullet \delta l)$. Then $\Phi(a^* l) = (\delta a^* \bullet \delta l) = \delta a$.
Since $\delta$ is a bijection, $\Phi$ is a surjection, and therefore $\Phi : (S, l) \cong (T, k)$. This completes the proof of Theorem 22.  ∎

**Corollary 23**. *Let [S, l] and [T, k] be laterally commutative semiheaps with bi-unitary elements l and k respectively. Then [S, l] $\cong$ [T, k] if and only if ([S, l]) $\cong$ ([T, k]).*

Proof. This follows from Theorem 22 and the fact that the mappings $\Phi$ and $\Gamma$ of Theorem 3 are mutually inverse mappings.  ∎

The proofs of the following Theorem and Proposition are straightforward and are omitted.

**Theorem 24** *If $\delta$ : [S] $\cong$ [T] and* **[S]** *is a semiheap* [*bi-unitary semiheap; a generalized heap; a laterally commutative semiheap; laterally commutative*] *then* **[T]** *is a semiheap* [*bi-unitary semiheap; a generalized heap; a laterally commutative semiheap; laterally commutative*]. *If S is right modular then so is T.*

**Proposition 25.** *A laterally commutative semiheap S satisfies the generalized heap axioms.*

Proof. [xx[yyz]] = [[yyz]xx] = [yy[zxx]] = [yy[xxz]] and [[xyy]zz] = [zz[xyy]] = [[zzx]yy] = [zz[xyy]] = [zz[yyx]]  ∎

## 5. Generalised heaps associated with groupoids determined by involutive, idempotent-fixed automorphisms on semilattices of groups.

**Definition 6.** An involutive homomorphism on a groupoid is a homomorphism whose square is the identity mapping. (Hence, an involutive homomorphism is an automorphism.)

**Definition 7:** A groupoid S is determined by a groupoid $(S, *)$ and an involutive automorphism $\alpha$ on $(S, *)$ if $xy = (\alpha x) * y$ $(x, y \in S)$.

**Definition 8.** A groupoid S is an ***inverse groupoid*** if for every $a \in S$ there exists a unique element $a^{-1} \in S$ such that $aa^{-1} \bullet a = a$ and $a^{-1} a \bullet a^{-1} = a^{-1}$. Then **[ ]**-1 is called the ***standard ternary operation on S***. and is denoted by **[ ]**-1 = **{ }.**

**Observation 1.** Note that since the idempotents of an inverse semigroup S commute, it is straightforward to prove that **[S]**-1 is a generalized heap.

**Theorem 26 ([10])** *A groupoid S satisfies*
(26.1) *S is an inverse groupoid,*
(26.2) $aa^{-1} = a^{-1}a \in E(S)$ *for all* $a \in S$,
(26.3) *there is an involutive automorphism* $\alpha$ *on S such that for all a, b, c* $\in S$, $ab \bullet c = (\alpha a) \bullet bc$, *and*
(26.4) *E(S) is a semilattice* **or** $(ab)^{-1} = \alpha b^{-1} \bullet \alpha a^{-1}$ *for all* $a, b \in S$
*if and only if S is determined by a semilattice of groups* $(S, *)$ *and an involutive, idempotent fixed automorphism* $\alpha$ *on* $(S, *)$. *In addition, in all such groupoids*
(26.5) $(ab \bullet c)d = a(bc \bullet d)$ *for all* $a,b,c,d \in S$.
(26.6) $\alpha a = a(aa^{-1})$ *for all* $a \in S$, *E(S) is a semilattice* **and** $(ab)^{-1} = \alpha b^{-1} \bullet \alpha a^{-1}$ *for all* $a,b \in S$.
(26.7) *if* $e \in E(S)$ *and* $a \in S$ *then* $ea = (\alpha a) \bullet e$.
(26.8) *for all* $a \in S$, *the inverse of a in* $(S, *)$ *is* $\alpha(a^{-1})$.

**Theorem 27.** *Let* $S(\alpha)$ *be a groupoid determined by a semilattice of groups* $(S, *)$ *and an involutive, idempotent-fixed automorphism* $\alpha$ *on* $(S, *)$. *Then the standard ternary operation on* $S(\alpha)$ *is a generalised heap. If* $S(\alpha)$ *has a left identity l then l is bi-unitary. Also, if* $S(\beta)$ *is a groupoid determined by the semilattice of groups* $(S, *)$ *and an involutive, idempotent-fixed automorphism* $\beta$ *on* $(S, *)$ *then* $\{S(\alpha)\} = \{S(\beta)\}$.

Proof. Note that since $ab = (\alpha a) * b$ $(a, b \in S)$, $\alpha(ab) = \alpha[(\alpha a) * b] = a * \alpha b = (\alpha a)(\alpha b)$, so $\alpha$ is an involutive isomorphism on $S(\alpha)$. Also, $a * b = (\alpha a) b$.

Then for any $a, b, c \in S$ we have $\{abc\} = (ab^{-1}) c = [(\alpha a) * b^{-1}] c = a * \alpha b^{-1} * c$. But, by (8) of Theorem 7, this product is equal to the standard ternary operation $\{abc\}$ on $(S, *)$, which is known to be a generalised heap [9]. Similarly, $\{abc\} = a * \beta b^{-1} * c$ in $S(\beta)$ and so, by (26.8), $\{S(\alpha)\} = \{S(\beta)\}$. ∎

**Corollary 28.** *Let* $S(\alpha)$ *be a groupoid determined by a semilattice of groups* $(S, *)$ *and an involutive, idempotent-fixed automorphism* $\alpha$ *on* $(S, *)$. *Then for any* $a, b, c \in S$, $\{abc\} = a * b^{-1} * c$, *where* $b^{-1}$ *is the inverse of b in* $(S, *)$.

**Theorem 29.** *Let* $S(\alpha)$ *and* $T(\beta)$ *be groupoids determined by the involutive, idempotent fixed automorphisms on the semilattice of groups* $(S, *)$ *and the semilattice of groups* $(T, \bullet)$ *respectively. Then* $\{S(\alpha)\} \cong \{T(\beta)\}$ *if and only if* $(S, *) \cong (T, \bullet)$.

Proof. ($\Leftarrow$) If $\delta: (S, *) \cong (T, \bullet)$ then by Corollary 28, for any $a, b, c \in S$, $\delta\{abc\} = \delta(ab^{-1}c) = \delta a(\delta b^{-1}) \delta c = \delta a (\delta b)^{-1} \delta c = \{(\delta a)(\delta b)(\delta c)\}$ and so $\{S(\alpha)\} \cong \{T(\beta)\}$.

($\Rightarrow$) Suppose that $\delta: \{S(\alpha)\} \cong \{T(\beta)\}$. For $a \in S$ we define $1_a$ as the identity of the group to which $a$ belongs in $(S, *)$. Recall that $1_{ab} = 1_a 1_b = 1_b 1_a$. Also, by Corollary 28, since for any $a,b,c \in S$, $\delta\{abc\} = \delta(ab^{-1}c) = \{(\delta a)(\delta b)(\delta c)\} = \delta a (\delta b)^{-1} \delta c$, the following equations are valid:
(1) $\delta(ab) = \delta a(\delta 1_a)^{-1} \delta b = \delta a(\delta 1_b)^{-1} \delta b$ and
(2) $\delta(a) = \delta a(\delta 1_a)^{-1} \delta 1_a = \delta 1_a (\delta a)^{-1} \delta 1_a$.

For $a \in S$ we define $\Phi(a) = \delta a (\delta 1_a)^{-1}$. We now prove that $\Phi: S(\alpha) \cong T(\beta)$. We will use without mention the facts that in a semilattice of groups idempotents are central and also $xx^{-1} = x^{-1}x$.

For $a,b \in S$ we have $\Phi(ab) = \delta(ab)(\delta 1_{ab})^{-1}$. From (1), $\Phi(ab) = \delta a(\delta 1_a)^{-1} \delta b (\delta 1_a 1_b)^{-1}$.
From (1), $(\delta 1_a 1_b)^{-1} = [\delta 1_a (\delta 1_a)^{-1} \delta 1_b]^{-1} = (\delta 1_b)^{-1} \delta 1_a (\delta 1_a)^{-1}$ and so
$\Phi(ab) = \delta a(\delta 1_a)^{-1} \delta b (\delta 1_a 1_b)^{-1} = \delta a(\delta 1_a)^{-1} \delta b (\delta 1_b)^{-1} \delta 1_a (\delta 1_a)^{-1}$. Therefore,
$\Phi(ab) = \delta 1_a (\delta 1_a)^{-1} \delta a(\delta 1_a)^{-1} \delta b (\delta 1_b)^{-1}$. Then, by (2),
$\Phi(ab) = \delta a(\delta 1_a)^{-1} \delta b (\delta 1_b)^{-1} = \Phi(a) \Phi(b)$ and $\Phi$ is a homomorphism.

If $\Phi(a) = \Phi(b)$ then $\delta a (\delta 1_a)^{-1} = \delta b (\delta 1_b)^{-1}$. By (2) and (1), this implies $\delta a = \delta a (\delta 1_a)^{-1} (\delta 1_a) =$ $= \delta b (\delta 1_b)^{-1} (\delta 1_a) = \delta (b1_a)$ and so $a = b1_a$. Similarly, $b = a1_b$. Therefore, $a = b1_a = a1_b 1_a = a1_a 1_b = a1_b = b$. Hence, $\Phi$ is one-to-one.

Since $\delta$ is surjective, for any $a \in S$ there exists $x \in S$ such that $\delta x = (\delta 1_a)^2$. Then by (2),
$$(3)\ \delta x = \delta x \delta 1_x (\delta 1_x)^{-1} = (\delta 1_a)^2 \delta 1_x (\delta 1_x)^{-1}.$$

Now $\Phi(1_a x) = \Phi(1_a) \Phi(x) = \delta 1_a (\delta 1_a)^{-1} \delta x (\delta 1_x)^{-1} = \delta 1_a (\delta 1_a)^{-1} (\delta 1_a)^2 (\delta 1_x)^{-1} = (\delta 1_a)^2 (\delta 1_x)^{-1} = \delta x (\delta 1_x)^{-1} = \Phi(x)$ and so $x = 1_a x$. Hence, $1_x = 1_a 1_x$. Also,
$$(4)\ \Phi(a1_x) = \Phi(a) \Phi(1_x) = \delta a (\delta 1_a)^{-1} \delta 1_x (\delta 1_x)^{-1}.$$

But, using (2) and (3), $(\delta 1_a)^{-1} = (\delta 1_a)^{-3} (\delta 1_a)^2 = (\delta 1_a)^{-3} \delta x = (\delta 1_a)^{-3} \delta x \delta 1_x (\delta 1_x)^{-1} = (\delta 1_a)^{-1} \delta 1_x (\delta 1_x)^{-1}$. Therefore, using (4), $\Phi(a1_x) = \delta a (\delta 1_a)^{-1} \delta 1_x (\delta 1_x)^{-1} = \delta a (\delta 1_a)^{-1} = \Phi(a)$.
So, $a = a1_x$ and $1_a = 1_a 1_x = 1_x$. Then, $\Phi(ax) = \Phi(a) \Phi(x) = \delta a (\delta 1_a)^{-1} \delta x (\delta 1_x)^{-1} = \delta a (\delta 1_a)^{-1} (\delta 1_a)^2 (\delta 1_a)^{-1} = \delta a (\delta 1_a)^{-1} (\delta 1_a) = \delta a$. Since $\delta$ is surjective, $\Phi$ is surjective. Hence, $\Phi : (S, *) \cong (T, \bullet)$. ∎

## 6. Idempotent automorphisms on right modular groupoids with a left identity element.

The proof of the following proposition is straightforward and is omitted.

**Proposition 30.** *Let $S$ be a right modular groupoid with left identity element $l$. Then the unary mapping $x \mapsto xl = x'$ is an involutive automorphism on $S$ and $S$ is $'$-unary.*

**Proposition 31.** *Let $S$ be a right modular groupoid with left identity element $l$. Then*

(31.1) $[S]'$ *is a laterally commutative and left commutative $l$-unital semiheap*;
(31.2) $(ab \bullet c)d = a(bc \bullet d)$ *for all $a,b,c,d \in S$* ;
(31.3) $l$ *is a lateral unit with respect to $[\ ]*$ if and only if $* = '$* ;
(31.4) $l$ *is bi-unitary with respect to $[\ ]*$ if and only if $l = l^*$* ;
(31.5) $k$ *is bi-unitary with respect to $[\ ]*$ implies $l = kk^* = k^*k$ and*
(31.6) *if a unary operation $*$ on $S$ is an involution then for all $a \in S$, $a = al$.*

Proof. Since $l' = ll = l$ it is clear that $l$ is a bi-unitary element. Also, $[lal]' = l(al) \bullet l = a$, so $l$ is a unit. Since $S$ is a right modular groupoid, $[\ ]'$ is laterally commutative. As noted in the paragraph above Proposition 1, $S$ is left lateral. It follows that $[S]'$ is left commutative. Then for all $a,b,c,d,e \in S$, $[[abc]' de]' =$
$= \{(ab' \bullet c)d'\}e = (d' c \bullet ab')e = \{a(d' c \bullet bl)\}e = a\{b(d' c \bullet l)\} \bullet e = a\{b(c \bullet d')\} \bullet e = a\{b(dc')\} \bullet e = a\{(dc' \bullet b)l\} \bullet e = [a [dcb]' e]'$. Then $[ab [cde]']' = [[cde]' ba]' = [[edc]' ba]' = [e [bcd]' a]' = [e[dcb]' a]' = [a [dcb]' e]'$. So we have proved that $[[abc]' de]' = [a [dcb]' e]' = [ab [cde]']'$; therefore $[S]'$ is a semiheap and (31.1) is valid.

Let $a,b,c \in S$. Then $(ab \bullet c)d = dc \bullet ab = ba \bullet cd = bc \bullet ad = (la)(bc \bullet d) = a (bc \bullet d)$, proving (31.2).

Since for all $a \in S$, $a = al \bullet l = [lal]'$, $l$ is a lateral unit with respect to $[\ ]'$. Conversely, if $l$ is a lateral unit with respect to $[\ ]*$ then $a = la^* \bullet l$ and so $a' = al = [la^* \bullet l]l = a^*$, which proves (31.3).

If $l$ is bi-unitary with respect to $[\ ]*$ then for all $a \in S$, $a = [lla]* = ll^* \bullet a = l^* a$ and, since left inverses are unique, $l = l^*$. Conversely, if $l = l^*$ then for all $a \in S$, $a = ll^* \bullet a = [lla]* = al^* \bullet l = [all]*$, proving (31.4).

If $k$ is bi-unitary with respect to $[\ ]*$ then for all $a \in S$, $a = [kka]* = kk^* \bullet a$, so $l = k^*k = kk^*$, proving (31.5).
If a unary operation $*$ on $S$ is an involution then, by definition, for all $a,b \in S$, $(a^*)^* = a$ and $(ab)^* = b^*a^*$.
Therefore, $l^* = (ll)^* = l^*l^* = (l\,l^*)\,l^* = (l^*l^*)l = (l^*l)$. So $l = (l^*)^* = (l^*l)^* = l^*\,(l^*)^* = l^*l = l^*$. Then
$al = (a^*)^*\,l^* = (la^*)^* = (a^*)^* = a$. This proves (31.6).  ∎

**Theorem 32.** *Let $(S, l)$ and $(T, k)$ be right modular groupoids with left identities $l$ and $k$ respectively. Then $[S] \cong [T]\,'$ if and only if there is a bijection $\delta : S \to T$ such that $\delta(ab) = \{\delta l\,(\delta a \bullet k)\}\,\delta b$, for all $a,b \in S$.*

Proof. ($\Rightarrow$) If $\delta : [S] \cong [T]\,'$ then for all $a,b,c \in S$, $\delta(ab \bullet c) = \delta[abc] = [(\delta a)(\delta b)(\delta c)]\,' = \{\delta a\,(\delta b \bullet k)\}\,\delta c$
Therefore, $\delta(ab) = \delta(la \bullet b) = \{\delta l\,(\delta a \bullet k)\}\,\delta b$.

($\Leftarrow$) Suppose that there is a bijection $\delta : S \to T$ such that $\delta(ab) = \{\delta l\,(\delta a \bullet k)\}\,\delta b$, for all $a,b \in S$. Then
$\delta(b) = \delta(lb) = \{\delta l\,(\delta l \bullet k)\}\,\delta b$ and, since left inverses are unique in a right modular groupoid, $k = \delta l\,(\delta l \bullet k)$.

Then, for all $a,b,c \in S$, $\delta[abc] = \delta(ab \bullet c) = \{\delta l\,(\delta(ab) \bullet k)\}\,\delta c$.
Now $\delta(ab) \bullet k = \{(\delta l\,(\delta a \bullet k))\,\delta b\}k = (k \bullet \delta b)\,(\delta l\,(\delta a \bullet k)) = \delta b(\delta l\,(\delta a \bullet k)) = \delta b(\delta a\,(\delta l \bullet k))$ and so
$\delta[abc] = \delta(ab \bullet c) = \{\delta l\,(\delta(ab) \bullet k)\}\,\delta = [\delta l\{\delta b(\delta a\,(\delta l \bullet k))\}]\,\delta c = [\delta b\{\delta l(\delta a\,(\delta l \bullet k))\}]\,\delta c =$
$= [\delta b\{\delta a\,(\delta l\,(\delta l \bullet k))\}]\,\delta c = \{\delta b\,(\delta a \bullet k)\}\,\delta c = \{\delta a\,(\delta b \bullet k)\}\,\delta c = [(\delta a)(\delta b)(\delta c)]\,'$. Hence,
$\delta : [S] \cong [T]\,'$.  ∎

**Corollary 33.** *Let $(S, l)$ and $(T, k)$ be right modular groupoids with left identities $l$ and $k$ respectively. If $[S] \cong [T]\,'$ then $al = a$, for all $a \in S$, and $S$ is a commutative monoid.*

Proof. Using Theorem 32, for all $a,b \in S$, $\delta(b) = \delta(lb) = \{\delta l\,(\delta l \bullet k)\}\,\delta b$ and so $k = \{\delta l\,(\delta l \bullet k)\}$. Then
$\delta(al) = \{\delta l\,(\delta a \bullet k)\}\,\delta l = \{\delta a\,(\delta l \bullet k)\}\,\delta l = \{\delta l\,(\delta l \bullet k)\}\,\delta a = \delta a$. Hence, $al = a$, for all $a \in S$. So $l$ is the identity of $S$. A right modular groupoid with identity is a commutative semigroup; so $S$ is a commutative monoid.  ∎

**Corollary 34.** *If $S$ is a right modular groupoid with left identitiy then $S$ is a commutative monoid if and only if $[S] \cong [T]\,'$ for some right modular groupoid $T$ with left identity.*

The properties of an $AG^*$ groupoid described in Note 1, Section 3, will be used in the proofs of Lemma 35 and Theorems 36, 37 and 38 below.

**Lemma 35.** *Let $S$ be an $AG^*$ groupoid. Then $S$ satisfies the identity $x^2 y^2 \bullet z = y^2 x^2 \bullet z$.*

Proof. Since right modularity implies mediality, $x^2 y^2 \bullet z = (xy \bullet xy)\,z = (yx \bullet yx)z = y^2 x^2 \bullet z$.  ∎

**Theorem 36.** *Let $S$ be an $AG^*$ groupoid. Define a ternary operation $[\ ]^*$ on $S$ as follows: for all $a,b,c \in S$, $[abc]^* = ac \bullet b$. Then $[S]^*$ is a semiheap that satisfies the identity $[abc]^* = [bac]^*$. Also, $[S]^*$ satisfies the generalised heap axioms.*

Proof. Let $a,b,c,d,e \in S$. Then $[[abc]^* de]^* = (ac \bullet b)e \bullet d = (de)(ac \bullet b) = (d \bullet ac)(eb) = (eb)(ac \bullet d) =$
$= (d \bullet eb)(ac) = (ac)(eb \bullet d) = (a \bullet eb)(cd) = (cd \bullet eb)a = (ce \bullet db)a = (a \bullet db)(ce) = (ce)(a \bullet db) = (a \bullet ce)(db) =$
$= (db)(ce \bullet a)$.

Also, $[ab\,[cde]^*]^* = a(ce \bullet d) \bullet b = (ce \bullet d)(ab) = (ce \bullet a)(db) = (db)(ce \bullet a)$.

Finally, $[a[dcb]^* e]^* = (ae)(db \bullet c) = (c \bullet ae)(db) = (db)(ae \bullet c) = (db)(ce \bullet a)$.

Hence, $[[abc]^* de]^* = [ab\,[cde]^*]^* = [a[dcb]^* e]^*$ and $[S]^*$ is a semiheap. Then, $[abc]^* = ac \bullet b = bc \bullet a = [bac]^*$ and so $[S]^*$ is a semiheap that satisfies the identity $[abc]^* = [bac]^*$.

Now $[[abb]^*cc]^* = (ab•b)c•c = c^2(b^2 a) = (ac^2)b^2 = b^2 c^2 •a$. Consequently, $[[acc]^*bb]^* = c^2 b^2 •a$ and, by Lemma 35, $[[abb]^*cc]^* = [[acc]^*bb]^*$.

Since, $[aa[bbc]^*]^* = a(bc•b)•a = (bc•b)a^2 = (cb^2)a^2 = a^2 b^2 •c$, $[bb[aac]^*]^* = b^2 a^2 •c$ and, by Lemma 35, $[aa[bbc]^*]^* = [bb[aac]^*]^*$. So we have proved that, $[S]^*$ satisfies the generalised heap axioms. ∎

**Theorem 37.** *Let S be an $AG^{**}$ groupoid or an $AG^{**}$ groupoid. Then* **[S]** *is a laterally commutative semiheap and satisfies the generalized heap axioms.*

Proof. This follows from Theorem 18 and Proposition 25. ∎

**Definition 9.** Let $\Pi$ be any permutation of the set $\{1,2,3\}$. Let [ ] be any ternary operation. We define $[\ ]_\Pi$ as follows: $[x_1 x_2 x_3]_\Pi = [x_{\Pi(1)} x_{\Pi(2)} x_{\Pi(3)}]$. We call the **[ ]$_\Pi$** *natural ternary operation induced by* $\Pi$. Note that **[ ]** = **[ ]**$_{(1,2,3)}$, where $\Pi = (1,2,3)$ is the identity permutation. We call $[\ ]_\Pi$ the ***dual*** of the ternary operation [ ] when $\Pi = (3,2,1)$. *We state without proof that a ternary operation is a semiheap [a generalised heap] if and only if its dual is a semiheap [generalised heap].*

**Theorem 38.** *Let S be an $AG^*$-groupoid. Then any ternary product on S defined by* $[a_1 a_2 a_3]_\Pi = a_{\Pi(1)} a_{\Pi(2)} •a_{\Pi(3)}$ *is a semiheap and satisfies the generalised heap axioms*, where $\Pi$ is any permutation of the set $\{1, 2, 3\}$.

Proof. As a result of Theorems 36 and 37, theorem 38 is valid for $\Pi = (1,2,3)$, $\Pi = (3,2,1)$, $\Pi = (1,3,2)$ and $\Pi = (2,3,1)$. We need only prove Theorem 38 for $\Pi = (2,1,3)$, since $[abc] = ba•c = ca•b$.

Now $[[abc]de] = d(ba•c)•e = (ba•dc)e$, $[ab[cde]] = ba(dc•e) = (dc•ba)e$ and $[a[dcb]e] = (cd•b)a•e = (ab•cd)e$. Using facts from Note 1, $[[abc]de] = [ab[cde]] = [a[dcb]e]$ and $[S]$ is a semiheap.

Using Lemma 35, it is straightforward to prove that $[S]$ satisfies the generalized heap axioms. ∎

**Proposition 39.** *If S is a right modular groupoid with left identity element l then Theorem 38 is <u>not</u> necessarily valid for S.*

Proof: In fact, we give an example of a right modular groupoid S with a left identity element *l*, where $[\ ]_{(1,3,2)}$ and $[\ ]_{(2,1,3)}$ are <u>not</u> semiheaps. Consider the groupoid *S* with the following Cayley table:

| S | a | b | c | l |
|---|---|---|---|---|
| a | l | c | b | a |
| b | b | l | a | c |
| c | c | a | l | b |
| l | a | b | c | l |

Then $S = \{a,b,c,l\}$ is a right modular groupoid with left identity element *l* (cf. [3], page 146). However, $[[abb]_{(1,3,2)} lb]_{(1,3,2)} = (ab•b)b•l = ab•l = ba = b \neq c = ab = al•b = [abl]_{(1,3,2)} = [ab [blb]_{(1,3,2)}]_{(1,3,2)}$. Also, $[[acb]_{(2,1,3)} aa]_{(2,1,3)} = a(ca•b)•a = a \neq l = (ca•c) = ca•(ab•a) = [ac [baa]_{(2,1,3)}]_{(2,1,3)}$. ∎

Note also that in this example *S* satisfies the equation $x^2 = l$. Therefore $[x[xyx]x] = x(xyx)•x = (xy•x^2)x =$

$= (xy \cdot l)x = (yx)x = x^2 y = ly = y$ and so [S] is outer lateral. Since $[xxx] = x^2 x = x$, by Theorem 37, [S] is a heap. While the left identity of a right modular groupoid is outer lateral, this example shows that an outer lateral element is not a left identity element in general.

In addition, if T is the dual groupoid of S then [T] is **_not_** a semiheap, since $[[blb]\, bl] = c \neq b = [b[bbl]\, l]$.

**Proposition 40.** *Let S be a right modular groupoid with left identity element l. Then*

(40.1) *S is '-unary, where $x \mapsto xl = x\,'$;*
(40.2) *for any unary operation $x \mapsto x\,'$ on S that satisfies the identities $x\,'\,' = x$ and $(xy)\,' = x\,'y\,'$,*
   $x\,'y = y\,'x$ *if and only if* $(xy)\,' = yx$;
(40.3) *the only unary operation $x \mapsto x\,'$ on S that satisfies the identities $x\,'\,' = x$, $(xy)\,' = x\,'y\,'$ and*
   $x\,'y = y\,'x$ *is* $x\,' = xl$ *and*
(40.4) *if $\Psi : S \to S$ is an automorphism satisfying $\Psi x \cdot y = \Psi y \cdot x$ for all $x,y \in S$ then $\Psi x = xl$.*

Proof: Proposition 30 is (40.1).
Now if $x \mapsto x\,'$ satisfies the identities $x\,'\,' = x$ and $(xy)\,' = x\,'y\,'$ then $x\,'y = y\,'x$ implies $yx = y\,'\,'x = $
$= x\,'y\,' = (xy)\,'$. Conversely, if $(xy)\,' = yx$ then $x\,'y = x\,' \cdot y\,'\,' = (x \cdot y\,')\,' = y\,'x$. This proves (40.2).
Now suppose that $x \mapsto x\,'$ satisfies the identities $x\,'\,' = x$, $(xy)\,' = x\,'y\,'$ and $x\,'y = y\,'x$. Then
$l\,' = (ll)\,' = l\,'l\,'$. So $l = (l\,')\,' = (ll\,')\,' = l\,'l = (l\,'l\,')l = ll\,' \cdot l\,' = l\,'l\,' = l\,'$. Then we have
$xl = xl\,' = lx\,' = x\,'$, which proves (40.3).

Finally, suppose that $\Psi : S \to S$ is an automorphism satisfying $\Psi x \cdot y = \Psi y \cdot x$ for all $x,y \in S$. Then
$\Psi l = \Psi (ll) = \Psi l \cdot \Psi l$. Since $\Psi x = \Psi (ll \cdot x) = (\Psi l \cdot \Psi l) \cdot \Psi x$, $\Psi l = \Psi l \cdot \Psi l = l$. Then
$\Psi x \cdot l = \Psi l \cdot x = x$ and so $xl = (\Psi x \cdot l)l = l \Psi x = \Psi x$, which proves (40.4). ∎

**Theorem 41.** *Let S be a groupoid with left identity element l. If S is '-unary then $x\,' = xl$. Furthermore, $[S]\,'$ is a semiheap if and only if S is a right modular groupoid.*

Proof. If S is '-unary then $l\,' = (ll)\,' = l\,'l\,' = (l\,')\,'l = l$.
Then $xl = xl\,' = lx\,' = x\,'$. Hence, $x\,' = xl$ and $x = (x\,')\,' = xl \cdot l = x\,'l$. This implies that $[lxl]\,' = lx\,' \cdot l = x$ and so $l$ is a lateral unit of $[S]\,'$.

Note also that $(xy)\,' = x\,'y\,' = (y\,')\,'x = yx$. Hence, $(xy)(zw) = \{(zw)(xy)\}\,' = (zw)\,'(xy)\,' = (wz)(yx)$.
Suppose that $[S]\,'$ is a semiheap. Then $[[lxl]\,'yz]\,' = [lx\,[lyz]\,']\,' = [l[ylx]\,'z]\,'$.

Now, since $l$ is a lateral unit, $[[lxl]\,'yz]\,' = [xyz]\,' = xy\,' \cdot z$. Then, $[lx\,[lyz]\,']\,' = [lx(y\,'z)]\,' = x\,' \cdot y\,'z$.
So $xy\,' \cdot z = x\,' \cdot y\,'z$ and, for any $w \in S$, setting $y = w\,'$, we have $xw \cdot z = x\,' \cdot wz = lx\,' \cdot wz$. But, using the fact that $(xy)(zw) = (wz)(yx)$, we have $xw \cdot z = lx\,' \cdot wz = zw \cdot x\,'l = zw \cdot x$, proving that S is right modular.

Conversely, if S is right modular then, by (31.1), $[S]\,'$ is a semiheap. ∎

**Corollary 42.** *Let S be a '-unary groupoid with left identity element l. Then S is right modular if and only if it is medial.*

Proof. ($\Rightarrow$) A right modular groupoid is medial.
($\Leftarrow$) By Theorem 41 and its proof, $x\,' = xl$, $x = (x\,')\,' = xl \cdot l = x\,'l$ and $(xy)\,' = x\,'y\,' = (y\,')\,'x = yx$.
So $xy \cdot z = (z \cdot xy)\,' = z\,' \cdot (xy)\,' = zl \cdot yx = zy \cdot lx = zy \cdot x$ and so S is right modular. ∎

**Theorem 43.** *Let S be a groupoid with left identity element l. If $[S]\,'$ is an l-unitary semiheap*

*then S is ʹ-unary.*

Proof. Firstly, $l\,ʹ = ll = l$. Then $x = [lxl]\,ʹ = xl \bullet l$, for all $x \in S$. That is, $(x\,ʹ)\,ʹ = x$. For any $x, y \in S$ we have $[lx\,[lyl]\,ʹ\,]\,ʹ = [lxy]\,ʹ = xl \bullet y = [l\,[ylx]\,l]\,ʹ = [l(yl\,ʹ \bullet x)\,l]\,ʹ = \{(yl\,ʹ \bullet x)l\}l = yl\,ʹ \bullet x = yl \bullet x$. That is, $x\,ʹ y = y\,ʹ x$.

Also, $xy = (xl \bullet l)y = (xl)\,ʹ y = y\,ʹ \bullet xl = y\,ʹ x\,ʹ = (yl)(xl)$. Finally, $[yl\,[(xl)ll]\,ʹ\,]\,ʹ = [yl(xl)]\,ʹ = yl\,ʹ \bullet xl = yl \bullet xl = [y\,[l(xl)l]\,ʹ l]\,ʹ = [y(xl)\,l]\,ʹ = y(xl \bullet l) \bullet l = (yx)l$. That is, $y\,ʹ x\,ʹ = (yx)\,ʹ = (x\,ʹ)\,ʹ y = xy$. ∎

**Corollary 44.** *Let S be a groupoid with left identity element l. Then the following are equivalent*: (1): $[S]\,ʹ$ *is an l-unitary semiheap* (2): *S is right modular;* (3): $[S]\,ʹ$ *is laterally commutative and* $(x\,ʹ)\,ʹ = x$ ($x \in S$); (4): $[S]$ *is an l-bi-unital semiheap and* (5): *S is paramedial. Consequently*, (6): *a groupoid T is a right modular groupoid with left identity element if and only if* $[T]$ *is an l-bi-unital semiheap for some idempotent element* $l \in T$. *Also*, (7): *S is ʹ-unary if and only if it is reversible.*

Proof. (1⇒2) By Theorem 43, ʹ satisfies the identities $(x\,ʹ)\,ʹ = x$, $(xy)\,ʹ = x\,ʹ y\,ʹ$ and $x\,ʹ y = y\,ʹ x$. Then by Theorem 41, *S* is right modular. (2⇒1) This follows from (31.1). (2⇒3) This follows from (31.1) and (40.1).

(3⇒2) For $x, y, z \in S$, we have $xy \bullet z = x(yl \bullet l) \bullet z = [x(yl)\,z]\,ʹ = [z(yl)\,x]\,ʹ = z(yl \bullet l) \bullet x = zy \bullet x$. Hence, (1), (2) and (3) are equivalent statements.

Assume (2) is valid. By Theorem 20, $[S]$ is an *l*-bi-unital semiheap in a right modular groupoid with left identity element *l*. So (2⇒4). Assume (4) is valid. We will prove that *S* is ʹ-unary. Let $x,y,z,w \in S$. Then $x = [xll] = xl \bullet l$; that is, $(x\,ʹ)\,ʹ = x$. Then, $(xy)l = [l\,[lxy]\,l] = [[lyx]\,ll] = [lyx] = yx$, so $xl \bullet y = [xl\,[yll]] = [x[lyl]\,l] = \{x \bullet yl\}l = yl \bullet x$; that is, $x\,ʹ y = y\,ʹ x$. Finally, $xl \bullet yl = [xl\,[lyl]] = [[xll]\,yl] = [xyl] = xy \bullet l$; that is, $(xy)\,ʹ = x\,ʹ y\,ʹ$. So we have proved that *S* is ʹ-unary. Therefore, as mentioned in the remarks following Definition 3, *S* is reversible. So then $xy \bullet zw = wz \bullet yx = [wz\,[xyl]] = [w[yxz]\,l] = \{w(yx \bullet z)\}l = (yx \bullet z)w$ and, setting $y = l$ in $wz \bullet yx = (yx \bullet z)w$ gives $wz \bullet x = xz \bullet w$; that is, *S* is right modular. So (4⇒2) and we have proved that (1), (2), (3) and (4) are equivalent statements. The equivalence of (2) and (5) is Proposition 1. We have already seen that a right modular groupoid is ʹ-unary and, therefore, is reversible.

Then, (6) follows from Theorem 20 and the equivalence of (2) and (4). Finally, we have seen in the comments following Definition 3 that if *S* is ʹ-unary then it is reversible. It is straightforward to prove that a reversible groupoid with a left identity is ʹ-unary. This proves the validity of (7). ∎

We know from Theorem 22 that two right modular groupoids with left identities are isomorphic if and only if they have isomorphic natural ternary operations. We now prove that the natural ternary product of a right modular groupoid *S* with left identity is isomorphic to the ternary operation of a second, arbitrary groupoid if and only if that groupoid is isomorphic to *S*. We need a lemma, the proof of which is straightforward and omitted.

**Lemma 45.** $[S]$ *is laterally commutative if and only if S is right modular.*

**Theorem 46.** *Suppose that S is a right modular groupoid with left identity l. Then* $S \cong T$ *if and only if* $[S] \cong [T]$.

Proof. (⇒) This follows from Proposition 3.
(⇐) By Theorem 20, $[S]$ is laterally commutative. By Theorem 24, $[T]$ is laterally commutative and, by Lemma 45, *T* is right modular. If $\delta : [S] \cong [T]$ then it is straightforward to prove that $(\delta\,l)^2$ is a left identity in *T*. It then follows from Theorem 22 that $S \cong T$. ∎

**Theorem 47.** *If S is a right modular groupoid with left identity l and* $[S]\,ʹ \cong [T]$ *then T is a commutative monoid.*

Proof. ($\Rightarrow$) By Proposition 31, [S]' is laterally commutative, with unit element *l*. Then, by Theorem 24, [T] is laterally commutative and, by Lemma 45, *T* is right modular. Since *l* is a unit of [S]' and $\delta : [S]' \cong [T]$, it follows that $(\delta l)^2$ is the left identity of *T*. By Corollary 33, *T* is a commutative monoid. ∎

## 7. Units in ternary operations of groupoids

**Theorem 48.** *Let S be a right modular groupoid with left identity element l. Then m is a lateral unit of* [S] *if and only if* $ml \bullet m = l$ *and* $x = xl$ ( $x \in S$ ).

Proof. ($\Rightarrow$) $[mxm] = x = mx \bullet m$ ( $x \in S$ ). So $ml \bullet m = l$. Then $x = mx \bullet m = mx \bullet lm = ml \bullet xm = x(ml \bullet m) = xl$.
($\Leftarrow$) For any $x \in S$, $mx \bullet m = mx \bullet lm = ml \bullet xm = x(ml \bullet m) = xl = x$. ∎

**Corollary 49.** *Let S be a right modular groupoid with left identity element l. Then m is a lateral unit of* [S] *if and only if* $l = m^2$ *and* $x = xl$ ( $x \in S$ ).

Proof. ($\Rightarrow$) By Theorem 48, $m = ml$ and $l = ml \bullet m = m^2$.
($\Leftarrow$) If $l = m^2$ and $x = xl$ then, for any $x \in S$, $x = xl = l(xl) = m^2 \bullet (xl) = (mx)(ml) = mx \bullet m$. ∎

**Corollary 50.** *Let S be a right modular groupoid with left identity element l. Then l is the only idempotent lateral unit of* [S].

Proof. If $m = m^2$ is a lateral unit of [S] then $m = m^2 = m^2 l = ml$. So $l = ml \bullet m = m^2 = m$. ∎

**Corollary 51.** *Let S be a right modular groupoid with left identity element l. The set LU of lateral units of* [S] *is a right modular subgroupoid of S, with left identity l or it is empty. If* $LU \neq \emptyset$ *then LU is a commutative monoid if and only if* $xl \bullet x = l$ *implies* $x^2 = l$ ( $x \in S$ ).

Proof. If $LU \neq \emptyset$, then let $m, n \in LU$. Then, by Theorem 48, $(mn)l \bullet (mn) = (ml \bullet nl)(mn) = (ml \bullet m)(nl \bullet n) = l$. Using Theorem 48 again, $mn \in LU$ and $l \in LU$. ∎

**Theorem 52.** *Let S be a right modular groupoid with left identity element l. Then m is a lateral unit of* [S]' *if and only if* $(ml \bullet m)l = l$ *if and only if* $ml \bullet m = l$ *if and only if* $m = mm^2$ *and* $l = m^2 m^2$.

Proof. For any $x \in S$, $[mxm]' = m(xl) \bullet m = ml \bullet (xl)m = (xl) \bullet (ml \bullet m) = \{(ml \bullet m)l\}x$. Hence, *m* is a lateral unit of [S]' if and only if $[mxm]' = x$ ( $x \in S$ ) if and only if $(ml \bullet m)l = l$. Now $(ml \bullet m)l = l$ implies $l = lm \bullet ml = m \bullet ml = ml \bullet m$. So $(ml \bullet m)l = l$ if and only if $ml \bullet m = l$. If $ml \bullet m = l$ then $m = lm = (ml \bullet m)m = m^2 \bullet ml = m^2 l \bullet ml = (m^2 m)l = mm^2$. Also, $ml = (mm^2)l = m^2 m$ and $l = ml \bullet m = (m^2 m)m = m^2 m^2$. Conversely, if $m = mm^2$ and $l = m^2 m^2$ then $ml \bullet m = (mm^2)l \bullet m = (m^2 m)m = m^2 m^2 = l$. ∎

The proofs of the Corollaries below are straightforward and are omitted.

**Corollary 53.** *Let S be a right modular groupoid with left identity element l. Then l is the only idempotent lateral unit of* [S]'.

**Corollary 54.** *Let S be a right modular groupoid with left identity element l. Then the set* (LU)' *of lateral units of* [S]' *is a right modular subgroupoid of S, with left identity element l. Also* (LU)' *is a commutative monoid if and only if* $ml \bullet m = l$ *implies* $m^2 = l$.

**Corollary 55.** *Let S be a right modular groupoid with left identity element l. Then a lateral unit of* [S]' *is a unit.*

**Theorem 56.** *Let S be a right modular groupoid with left identity element l. Then m is an outer lateral unit of* [S] *if and only if* $l = m^2 m^2$.

Proof. For any $x \in S$, $[m[mxm]\,m] = m(mx \bullet m) \bullet m = (ml)\{(mx \bullet m)m\} = (ml)(m^2 \bullet mx) = (mm^2)(mx) =$
$= m^2 \bullet m^2 x = m^2 l \bullet m^2 x = (m^2 m^2)\,x$ and so $x = [m[mxm]\,m]$ ($x \in S$) if and only if $l = m^2 m^2$. ∎

**Corollary 57.** *Let S be a right modular groupoid with left identity element l. Then l is the only idempotent outer lateral unit of* $[S]$.

**Corollary 58.** *Let S be a right modular groupoid with left identity element l. The set OLU of outer lateral units of* $[S]$ *is a right modular subgroupoid of S, with left identity element l. Also, OLU is a commutative monoid if and only if* $l = m^2 m^2$ *implies* $m = m^2 (mm^2)$.

**Theorem 59.** *Let S be a right modular groupoid with left identity element l. Then m is an outer lateral unit of* $[S]\,'$ *if and only if* $l = (ml \bullet m)(m \bullet ml)$.

Proof. For any $x \in S$, $[m\,[mxm]\,'\,m]\,' = [m\,\{m(xl) \bullet m\}\,m]\,' = m\{m \bullet m(xl)\} \bullet m = (ml)\,\{\{m \bullet m(xl)\} \bullet m\,\} =$
$= (ml) \bullet (ml)\{m(xl) \bullet m\} = (ml) \bullet (ml)\{(ml) \bullet (xl)\,m\} = (ml) \bullet (ml)\{(xl) \bullet (ml)\,m\} = (ml) \bullet (xl)\{(ml) \bullet (ml)\,m\} =$
$= (xl) \bullet (ml)\{(ml) \bullet (ml)\,m\} = \{(ml) \bullet (ml)\,m\}(ml) \bullet x = \{(ml \bullet m)(m \bullet ml)\}x$. So $x = [m\,[mxm]\,'\,m]\,'$ ($x \in S$) if and only if $l = (ml \bullet m)(m \bullet ml)$. ∎

**Corollary 60.** *Let S be a right modular groupoid with left identity element l. Then l is the only idempotent outer lateral unit of* $[S]\,'$.

**Corollary 61.** *Let S be a right modular groupoid with left identity element l. Then m is an outer lateral unit of* $[S]\,'$ *if and only if ml is an outer lateral unit of* $[S]\,'$.

**Corollary 62.** *Let S be a right modular groupoid with left identity element l. The set* $(OLU)\,'$ *of outer lateral units of* $[S]\,'$ *is a right modular subgroupoid of S, with left identity element l. Also,* $(OLU)\,'$ *is a commutative monoid if and only if* $l = (ml \bullet m)(m \bullet ml)$ *implies* $l = m^2 m^2$.

**Corollary 63.** *Let S be a right modular groupoid with left identity element l. Then* $OLU \cap (OLU)\,'$ *is a commutative submonoid of S, with identity element l.*

## 8. Characterisations of ternary operations on classes of groupoids

**Theorem 64.** $[S]$ *is the natural ternary product of a right modular groupoid with left identity element l if and only if* $[S]$ *is a laterally commutative, l-bi-unital semiheap.*

Proof. ($\Rightarrow$) This follows from Theorem 20.

($\Leftarrow$) The fact that the $\Phi$ and $\Gamma$ in the proof of Theorem 20 are mutually inverse mappings implies that $[S]$ is the ternary product of a right modular groupoid with left identity element $l$. ∎

**Theorem 65.** $\{S\} = [S]\,'$ *for some right modular groupoid S with left identity element l if and only if* $\{S\}$ *is a laterally commutative, left commutative and l-unital semiheap satisfying the identity* $\{xy\char`\^ l\}\char`\^ = \{x\char`\^ yl\}$ *for some idempotent unary operation* $\char`\^$.

Proof. ($\Rightarrow$) Suppose that $\{S\} = [S]\,'$ for some right modular groupoid $S$ with left identity element $l$. Then from (31.1), $\{S\}$ is a laterally commutative and left commutative $l$-unital semiheap.
Then, $\{xy\,'l\} = [xy\,'l]\,' = [x(yl)l]\,' = xy \bullet l = yx$. Then, since $a\,' = al$, $\{xy\,'l\}\,' = (yx)\,' = (yx)l =$
$= (xy)l \bullet l = (xl \bullet yl)l = [x\,'yl]\,' = \{x\,'yl\}$. Setting $' = \char`\^$ completes the proof of ($\Rightarrow$).

($\Leftarrow$) Define $a \bullet b = \{ba^\wedge l\}$. Now $l = (l^\wedge)^\wedge = \{ll^\wedge l\}^\wedge = \{l^\wedge ll\} = l^\wedge$. Thus, $l \bullet b = \{bl^\wedge l\} = \{bll\} = b$ and $l$ is a left identity element of $\{S, \bullet\}$ We now prove that $\{S, \bullet\}$ is right modular. Let $a, b, c \in S$. Then $ab \bullet c = \{c\{ba^\wedge\}^\wedge l\} = \{c\{b^\wedge al\}l\} = \{\{cla\}b^\wedge l\} = \{\{lca\}b^\wedge l\} = \{l\{b^\wedge ac\}l\} = \{b^\wedge ac\} = \{cab^\wedge\} = \{acb^\wedge\} = cb \bullet a$. We now prove that $\{abc\} = [S]^\wedge$. We have $[abc]^\wedge = ab^\wedge \bullet c = \{b^\wedge a^\wedge l\} \bullet l = \{c\{b^\wedge a^\wedge l\}^\wedge l\} = \{c\{bal\}l\}$. By Lemma 19, $[abc]^\wedge = \{c\{bal\}l\} = \{\{cll\}\{bal\}l\} = \{\{abc\}ll\} = \{abc\}$. This completes the proof of ($\Leftarrow$) and, therefore, of Theorem 65. ∎

**Theorem 66.** $\{S\} = [S] *$ *for some* $*$-*unary groupoid S with left identity element l if and only if* $\{S\}$ *is l-unital, left commutative, left l-consistent and satisfies the identities* $\{xly\} = \{ylx\}, l = l^\wedge$ *and* $\{xyz\}^\wedge = \{x^\wedge y^\wedge z^\wedge\}$ *for some idempotent unary operation* ^ *on S.*

Proof. ($\Rightarrow$) Suppose that $\{S\} = [S] *$ for some $*$-unary groupoid $S$ with left identity element $l$. Then $x^* l = (l x^*)^* = x$ and so $x^* = (x^* l)^* = xl^{**} = x(ll)^* = xl$. So $\{lxl\} = [lxl] * = x^* l = x = lx = [llx] * = \{llx\} = (x^*)^* = xl \bullet l = [xll] * = \{xll\}$. Therefore, $\{S\}$ is $l$-unital. Also, $\{xyz\} = [xyz] * = xy^* \bullet z = yx^* \bullet z = [yxz] * = \{yxz\}$, so $\{S\}$ is left commutative. Then $\{l\{lxy\}z\} = [l[lxy] * z] * = [l(x^* y)z] * = (x^* y)^* \bullet z = xy^* \bullet z = [xyz] * = \{xyz\}$ and so $\{S\}$ is left $l$-consistent. Also, $\{xly\} = [xly] * = xl \bullet y = x^* y = y^* x = yl \bullet x = [ylx] * = \{ylx\}$. Note also that $l^* = (ll)^* = l^* l^* = ll = l$. Finally, $\{xyz\}^* = ([xyz] *)^* = (xy^* \bullet z)^* = x^* y \bullet z^* = [x^* y^* z^*] * = \{x^* y^* z^*\}$. Setting $* = \wedge$ completes the proof of ($\Rightarrow$).

($\Leftarrow$) Define $xy = \{ylx^\wedge\}$. Since, $l^\wedge y = \{yll^{\wedge\wedge}\} = \{yll\} = y$, $l^\wedge = l$ is a left identity element. Then $(xy)^\wedge = \{ylx^\wedge\}^\wedge = \{y^\wedge l^\wedge x\} = \{y^\wedge lx\} = x^\wedge y^\wedge$. Also, $x^\wedge y = \{ylx\} = \{xly\} = y^\wedge x$ and so $S$ is $\wedge$-unary. So we need only show that $\{xyz\} = [xyz]^\wedge$; that is, that $\{xyz\} = xy^\wedge \bullet z$. Now $xy = \{ylx^\wedge\} = \{x^\wedge ly\} = y^\wedge x^\wedge$. So $xy^\wedge \bullet z = z^\wedge \bullet (xy^\wedge)^\wedge = z^\wedge \bullet (x^\wedge y) = \{(x^\wedge y)lz\} = \{\{ylx\}lz\} = \{l\{ylx\}z\} = \{l\{lyx\}z\} = \{yxz\} = \{xyz\}$. ∎

**Theorem 67.** $\{S\} = [S]$ *for some* $*$-*unary groupoid S with left identity element l if and only if* $\{S\}$ *is l-bi-unital, left l-consistent and* ^ -*congruent, l is a central commutant and* $\{S\}$ *satisfies the identity* $\{lxl\} = x^\wedge$ *for some unary operation* ^ *on S.*

Proof. ($\Rightarrow$) Suppose that $\{S\} = [S]$ for some $*$-unary groupoid $S$ with left identity $l$. Then $\{lxl\} = [lxl] = lx \bullet l = xl$. By Theorem 41, $\{lxl\} = xl = x^*$. Also, $\{llx\} = [llx] = ll \bullet x = x = (x^*)^* = xl \bullet l = [xll] = \{xll\}$ and so $\{S\}$ is $l$ bi-unital. In addition, $\{xly\} = [xly] = xl \bullet y = x^* y = y^* x = [ylx] = \{ylx\}$, so $l$ is a central commutant. We have $\{l\{lxy\}z\} = [l[lxy] z] = xy \bullet z = [xyz] = \{xyz\}$, so $\{S\}$ is left $l$-consistent. Finally, $\{xyz\}^* = [xyz]^* = (xy \bullet z)^* = (xy)^* z^* = (x^* y^*) \bullet z^* = [x^* y^* z^*] = \{x^* y^* z^*\}$, so $\{S\}$ is $*$-congruent. Setting $\wedge = *$ completes the proof of ($\Rightarrow$).

($\Leftarrow$) Suppose that $\{S\}$ is $l$ bi-unital, left $l$-consistent an $\wedge$-congruent, $l$ is a central commutant and $\{S\}$ satisfies the identity $\{lxl\} = x^\wedge$. Define $xy = \{lxy\}$. Since $\{S\}$ is $l$ bi-unital, $x = \{llx\} = lx$ and so $l$ is a left identity element. Also, $l = \{lll\} = l^\wedge$ and $xl = \{lxl\} = x^\wedge$. Then $(x^\wedge)^\wedge = (xl)l = \{l\{lxl\}l\}$ and, by left $l$-consistency, $(x^\wedge)^\wedge = \{l\{lxl\}l\} = \{xll\} = x$. Using left $l$-consistency again, $(xy)^\wedge = (xy)l = \{l\{lxy\}l\} = \{xyl\}$. Then $xy = \{xyl\}^\wedge = \{x^\wedge y^\wedge l^\wedge\} = \{x^\wedge y^\wedge l\} = \{l\{lx^\wedge y^\wedge\}l\} = (x^\wedge y^\wedge)l = (x^\wedge y^\wedge)^\wedge$ and so $(xy)^\wedge = (x^\wedge y^\wedge)$. Finally, $x^\wedge y = \{l\{lxl\}y\} = \{xly\}$. Since $l$ is a central commutant, $x^\wedge y = \{xly\} = \{ylx\} = y^\wedge x$. So $S$ is $*$-unary and $\{xyz\} = \{l\{lxy\}z\} = xy \bullet z = [xyz]$. This completes the proof of ($\Leftarrow$) and, therefore, of Theorem 67. ∎

We can use the proofs of Theorems 64, 65, 66 and 67 to produce results similar to Theorem 20. That is, there are bijections between the collection of right modular groupoids with a left identity element, the collection of laterally commutative, bi-unitary semiheaps and the collection of laterally commutative, left commutative and $l$-unital semiheaps satisfying the identity $\{xy^* l\}^* = \{x^* yl\}$ for some idempotent unary operation $*$.

Similarly, there are bijections between the collection of $*$-unary groupoids with a left identity element $l$, the collection of $l$-bi-unital, left $l$-consistent and $*$-congruent ternary operations, with $l$ a central commutant, satisfying the identity $\{lxl\} = x^*$ and the collection of $l$-unital, left commutative, left $l$-consistent ternary operations satisfying the identities $\{xly\} = \{ylx\}$ and $\{xyz\}^* = \{x^* y^* z^*\}$ for some idempotent unary operation $*$.

One can also attempt to characterise the ternary products of the groupoids dual to the classes of groupoids we have already investigated. For the duals of right modular groupoids with a left identity element we obtain the following result.

**Theorem 68.** *A ternary product { } is the natural ternary product of the dual groupoid* $(S,\times)$ *of a right modular groupoid* $(S,\bullet, l)$ *with left identity element l if and only if { } is right commutative, l-right unital and satisfies the identities* $\{x\{yzq\}w\} = \{zw\{yxq\}\}$ *and* $\{\{xyz\}qw\} = \{\{xyq\}wz\}$.

Proof. ($\Rightarrow$) First note that $l$ is a right identity of $(S,\times)$. Let $x, y, z, q, w \in S$. Then, by definition, $x \times y = y \bullet x$. Then $[xll] = (xl)l = l\bullet(l\bullet x) = x$, so $[(S,\times)]$ is $l$-right unital. Also, $[xyz] = (x\times y)\times z = z\bullet(y\bullet x) = y\bullet(z\bullet x) = $ $ = (x\times z)\times y = [xzy]$ and so $[(S,\times)]$ is right commutative. Since $(xy)\bullet(zw) = (w\times z)\times(y\times x) = (xz)\bullet(yw) = $ $ = (w\times y)\times(z\times x)$, $(S,\times)$ is also medial.

Consider that $[x[yzq]\,w] = x[(yz)q]\times w = w\bullet[q(zy)\bullet x] = q(zy)\bullet wx = qw\bullet(zy\bullet x) = qw\bullet(xy\bullet z) = (q\bullet xy)\bullet(wz) = $ $ = (z\times w)\times[(y\times x)\times q] = [zw\,[yxq]]$. Note that $(x\times y)\times z = z\bullet(y\bullet x) = y\bullet(z\bullet x) = (x\times z)\times y$. So then $[[xyz]\,qw] = [(x\times y)\times z](q\times l)\times w = [(x\times y)\times q](z\times l)\times w = [(x\times y)\times q]z\times w = [(x\times y)\times q]w\times z = [[xyq]\,wz]$. So $[(S,\times)]$ satisfies the two required identities. This completes the proof of ($\Rightarrow$).

($\Leftarrow$) Let $x, y, z, q, w \in S$. Define a product $\bullet$ on $S$ as follows: $x\bullet y = \{yxl\}$. Then $l\bullet y = \{yll\} = y$ and so $l$ is a left identity of $(S,\bullet)$. Also, $(x\bullet y)\bullet z = \{yxl\}\bullet z = \{z\{yxl\}l\} = \{xl\{yzl\}\} = \{x\{yzl\}l\} = (z\bullet y)\bullet x$ and so $(S,\bullet)$ is right modular, with left identity $l$.

We now prove that $\{\ \}$ is the natural ternary product of $(S,\times)$, the dual groupoid of $(S,\bullet)$. Now consider $[(S,\times)]$. We have $[xyz] = (x\times y)\times z = z\bullet(y\bullet x) = \{\{xyl\}zl\} = \{\{xyz\}ll\} = \{xyz\}$. This completes the proof of ($\Leftarrow$) and, therefore of Theorem 68. ∎

Theorems 46 and 47 lead us to conjecture that $[S] \cong [T]$ for some $*$-unary groupoid $S$ with a left identity element if and only if $S \cong T$ and that $[S]* \cong [T]$ for some $*$-unary groupoid $S$ with a left identity element only for a restricted class of groupoids $T$ (as in Theorem 47). We conjecture that there exists a right modular groupoid $S$ with left identity element and a commutative monoid $T$ such that $[S]' \cong [T]$ but $S$ and $T$ are **not** isomorphic. If this conjecture is true then it implies that $[S]' \cong [T]'$ does **not** imply that $S$ and $T$ are isomorphic.

If $S$ is an $AG^*$ or $AG^{**}$-groupoid without left identity element then it is difficult to see how to characterise $[S]$. This is because, without a special element, it is not clear how to define a product in terms of $[S]$. A lateral unit of $[S]$ is a right or left identity element and so $[S]$ has no such element. In an $AG^*$ or an $AG^{**}$-groupoid the square of the square of an outer lateral unit is a left identity and so such elements cannot be used to define a product using $[S]$. We prove these statements in the following Theorems.

**Theorem 69.** *If S is an $AG^*$ or $AG^{**}$-groupoid without left identity element then* $[S]$ *has no lateral or outer lateral units*.

Proof. Suppose that $S$ is an $AG^*$- groupoid without a left identity. Then the identity $xy\bullet z = y\bullet xz = y\bullet zx$ holds in $S$. If $l$ is a lateral unit then $x = [lxl] = lx\bullet l = xl^2$. Since $S$ is a right modular groupoid with right identity $l^2$, $S$ is commutative. Then $S$ has a (left) identity element $l^2$, a contradiction. If $l$ is an outer lateral unit then we have $x = [l[lxl]\,l] = l(lx\bullet l)\bullet l = (lx\bullet l)l^2 = (l^2\,l)(lx) = (l^2\,l)(xl) = (l^2\,x)l^2 = x(l^2)^2$. So, as before, $(l^2)^2$ is a right and, therefore, a left identity element, a contradiction.

If $S$ is an $AG^{**}$-groupoid (without left identity) then $S$ satisfies the identity $x\bullet yz = y\bullet xz$. If $l$ is a lateral unit then $S$ satisfies the identity $x = [lxl]$. Then $x = [lxl] = lx\bullet l = lx\bullet(l^2\,l) = l^2\,(lx\bullet l) = l^2\,x$, a contradiction. If $l$ is an outer lateral unit then we have $x = [l[lxl]\,l] = l(lx\bullet l)\bullet l = (lx\bullet l^2)l = (l^2\bullet xl)l = (x\bullet l^2\,l)l = \{l(l^2\,l)\}x = (l^2)^2\,x$, a contradiction. ∎

**Theorem 70.** *Suppose that S is a ∗-unary groupoid with a left identity element l. If S is left lateral then S is a right modular groupoid. If S is weakly associative then S is a commutative monoid.*

Proof. Now, since $S$ is ∗-unary, by Theorem 41, $x^* = xl$. Therefore, $x = (x^*)^* = (xl)l$ and $(xy)^* = (xy)l = x^* y^* = (y^*)^*(x^*)^* = yx$. If $S$ is left lateral then, by definition, it satisfies the identity $x \bullet y\, z = y \bullet xz$. Then $xy \bullet zw = (xy)\{(zw)l \bullet l\} = (zw)l \bullet (xy)l = wz \bullet yx$. It follows that $xy \bullet z = xy \bullet lz = zl \bullet yx = y\,(zl \bullet x) = y \bullet z^* x = y \bullet x^* z = zy \bullet x$ and so $S$ is right modular.

If $S$ is weakly associative then, by definition, it satisfies the identity $xy \bullet z = y \bullet xz$. Then $x = xl \bullet l = l \bullet xl = xl$. Also, $xy = (xy)l = y \bullet xl = yx$. Then $xy \bullet z = yx \bullet z = x \bullet yz$ and so $S$ is a commutative semigroup with identity element $l$.  ∎

The following Theorem follows straightforwardly from Theorem 41 and we omit the proof.

**Theorem 71.** *Let S be a ∗-unary groupoid with left identity l. If S satisfies any one of the identities $xy \bullet z = y^* \bullet xz$, $xy \bullet z = y \bullet xz^*$, $x \bullet y\, z = y^* \bullet xz$, $x \bullet yz = y \bullet x^* z$ or $x \bullet y\, z = y \bullet xz^*$ then S is a commutative monoid. If S satisfies the identity $xy \bullet z = y \bullet x^* z$ then S is right modular.*

**Example 2.** We conclude this section with an example of a ∗-unary groupoid $S$ with a left identity element $l$ that is **_not_** right modular. In this case $* = 1_S = {}^{\prime}$ and the Cayley table of $S$ is as follows:

| S | x | y | l |
|---|---|---|---|
| x | y | x | x |
| y | x | x | y |
| l | x | y | l |

## 9. Heaps, groups and Ward quasigroups

Vagner [13] proved that [S] is a heap if and only if every element is bi-unitary and [S] satisfies the identity $[[xyz]\,uv] = [xy\,[zuv]]$. We give another characterization of a heap, as follows:

**Theorem 72.** *[S] is a heap if and only if every element is bi-unitary and [S] satisfies*
*(1) $[xyz] = q$ implies $[zqx] = y$ and (2) $[xyz] = [[xyq]\,qz] = [x[qzy]\,q]$.*

Proof. ($\Rightarrow$) If [S] is a heap then, by definition, every element is bi-unitary and [S] satisfies the identity $[[xyz]\,uv] = [x[uzy]\,v] = [xy\,[zuv]]$. Setting $u = v = q$ and gives $[xyz] = [x[qzy]\,q]$. Also, $[xyz] = [[xy\,[qqz]] = [[xyq]\,qz]$. So [S] satisfies (2). In addition this implies that $[z[xyz]\,x] = [[zzy]\,xx] = y$. Therefore, $[xyz] = q$ implies $[zqx] = y$ and [S] satisfies (1). This completes the proof of necessity.

($\Leftarrow$) Assume that [S] satisfies (1) and (2). Let $[zuv] = q$, which implies $[vqz] = u$. So $[[xyz]\,uv] = [[xyz][vqz]\,v] = [[xyz]\,zq] = [xyq] = [xy\,[zuv]]$ and therefore [S] is a heap.  ∎

**Corollary 73.** *In a heap $[xyz] = x$ if and only if $y = z$ and $[xyz] = z$ if and only if $x = y$.*

Proof. By Theorem 72, a heap satisfies the implication $[xyz] = q$ implies $[zqx] = y$ and Corollary 72 follows readily from this fact.  ∎

**Proposition 74.** *In a heap [S], $[xyz] = [xwz]$ if and only if $y = w$, $[xyz] = [wyz]$ if and only if $x = w$ and $[xyz] = [xyw]$ if and only if $z = w$.*

Proof. If $[xyz] = [xwz]$ then $y = [[zzy]\,xx] = [z[xyz]\,x] = [z[xwz]\,x] = [[zzw]\,xx] = w$. If $[xyz] = [wyz]$ then $x = [xy\,[zzy]] = [[xyz]\,zy] = [[wyz]\,zy] = [wy\,[zzy]] = w$. If $[xyz] = [xyw]$ then $y = [[zzy]\,xx] = [z[xyz]\,x] = [z[xyw]\,x] = [[zwy]\,xx] = [zwy]$, which by Corollary 73 implies $z = w$.  ∎

We gather work from [2], [3], [6], [11] and [12] into the following results, thereby showing the close connection between groups, heaps and Ward quasigroups.

**Definition 9.** A *Ward groupoid* is a groupoid satisfying the identity $xz \bullet yz = xy$. A *Ward quasigroup* is a cancellative Ward groupoid. A groupoid is ***right solvable*** if for any two elements $a$ and $b$ there exists an element $x$ such that $a = bx$. $\mathbf{G}$ is the collection of all groups, $\mathbf{H}$ is the collection of all heaps and $\mathbf{WQ}$ is the collection of all Ward quasigroups.

**Proposition 75.** *The following statements are equivalent*:

(75.1) *S is a Ward quasigroup*;
(75.2) *S is a right solvable Ward groupoid and*
(75.3) *S is a quasigroup satisfying the identities* $x^2 = y^2$ *and* $xy \bullet z = x \bullet (z \bullet w^2 y)$.

Proof. Assume that (75.3) is valid. Let $e = x^2$ (any $x \in S$). The for any $x \in S$, $xe \bullet e = x \bullet (e \bullet ee) = x \bullet e$. Since $S$ is a quasigroup, this implies $xe = x$. For any $y, z \in S$, $yz \bullet ez = y \bullet (ez \bullet ez) = ye = y$. Hence, for any $x, y, z \in S$, $xz \bullet yz = x \bullet (yz \bullet ez) = xy$. So $S$ is a Ward groupoid and (75.1) is valid. Thus, $e \bullet ex = x^2 \bullet ex = x (ex \bullet ex) = xe = x$. Hence, $x = e \bullet ex = y^2 \bullet ex = y (ex \bullet ey)$, so $S$ is right solvable and we have proved that (75.2) is valid.

Assume (75.1) is valid. Then $x^2 = x^2 \, x^2$ and so $E(S) = \{x \in S : xx = x\}$ is non-empty. If $e, f \in S$ then $ef = ef \bullet ff = ef \bullet f$ which implies $ee = e = ef$ which implies $e = f$. So we can write $E(S) = \{e\}$. So $x^2 = y^2$ holds in $S$. Since $xe = xe \bullet ee = xe \bullet e$, $x = xe$. Then $xy \bullet z = (xy \bullet ey)(z \bullet ey) = xe \bullet (z \bullet ey) = x \bullet (z \bullet w^2 y)$. Hence, (75.3) is valid. So we need only prove that (75.2) implies (75.3).

Assume (75.2) is valid. For any $x, y, z \in S$ we have $x = yw$ for some $w$. Then, $x^2 = yw \bullet yw = y^2$. Also, if $e = x^2$ Then $xe = yw \bullet ww = yw = x = xx \bullet ex = e \bullet ex$. Also, $xy = xx \bullet yx = e(yx)$. Now if $ax \bullet bx$, $e = ax \bullet bx = bx \bullet ax = = ab = ba$. Then $a = e(ea) = ba \bullet ea = be = b$, so $S$ is right cancellative. If $xa = xb$ then $ax = e(xa) = e(xb) = bx$ and so $S$ is left cancellative. So $S$ is a quasigroup. Finally, $xy \bullet z = (xy \bullet ey) \bullet (z \bullet ey) = xe \bullet (z \bullet ey) = x \bullet (z \bullet w^2 y)$, so (75.2) implies (75.3). ∎

**Corollary 76.** $\mathbf{WQ}$ *is a variety of groupoids determined by the identities* $x^2 = y^2$, $xy = xz \bullet yz$ *and* $xy^2 = x$.

Proof. We write $x^2 = e$. Then $e \bullet xy = yy \bullet xy = yx$. Also, $x = xy^2 = xe = x^2 \bullet ex = e \bullet ex$. Then if $xy = xq$ then $yx = e \bullet xy = e \bullet xq = qx$. So $yq = yx \bullet qx = e$ and $e = ee = e \bullet yq = qy$. So $y = e \bullet ey = qy \bullet ey = qe = q$. If $yx = qx$ then $xy = e \bullet yx = e \bullet qx = xq$ and so $y = q$. Hence, we have a (Ward) quasigroup. Conversely, by Theorem 75 and its proof, every Ward Quasigroup satisfies the equations $x^2 = y^2$, $xy^2 = x$ and $xy = xz \bullet yz$. ∎

**Definition 10.** Define $\Psi : \mathbf{G} \to \mathbf{H}$ as follows: $\Psi(G, \bullet, e) = \{G\}$, where $\{xyz\} = x \bullet y^{-1} \bullet z$. Define $\Omega : \mathbf{H} \to \mathbf{G}$ as follows: $\Omega(H) = (H, \square_e)$, where $x \square_e y = [xey]$, $e \in H$. Define $\Pi : \mathbf{H} \to \mathbf{WQ}$ as follows: $\Pi([H]) = (H, \square_e)$, where $x \square_e y = [xye]$ and $e \in H$. Define $\Lambda : \mathbf{WQ} \to \mathbf{H}$ as follows: $\Lambda(W, \bullet, e) = \{W\}$, where $\{xyz\} = xy \bullet ez$. Define $\Gamma : \mathbf{WQ} \to \mathbf{G}$ as follows: $\Gamma(W, \bullet, e) = (W, \circ_e)$, where $x \circ_e y = x \bullet (e \bullet y)$. Define $\Phi : \mathbf{G} \to \mathbf{WQ}$ as follows: $\Phi(G, \bullet, e) = (G, *)$, where $x * y = x \bullet y^{-1}$.

**Note 2.** Without proof we observe that $\Psi$, $\Omega$, $\Pi$, $\Lambda$, $\Gamma$ and $\Phi$ are well defined: $(H, \square_e)$ is a group with identity $e$ and $x^{-1} = [exe]$, $(H, \square_e)$ is a Ward quasigroup with right unit $e$, $(W, \circ_e)$ is a group with identity element $e$ and $x^{-1} = e \bullet x$ and $(G, *)$ is a Ward quasigroup with right unit $e$.

**Proposition 77.** *For* $e, f \in H$, $\delta : (H, \square_e) \cong (H, \square_f)$, *where* $\delta x = [fex]$.

Proof. Now $\delta(x \square_e y) = \delta [xey] = [fe[xey]] = [fe[[xff] ey]] = [fe[xf[fey]]] = [[fex] f [fey]] = \delta x \square_f \delta y$ and so $\delta$ is a homomorphism. If $\delta x = \delta y$ then $[fex] = [fey]$ and so $x = [eex] = [[eff] ex] = [ef[fex]] =$
$= [ef[fey]] = [[eff] ey] = [eey] = y$ and so $\delta$ is an injection. Then $\delta [efx] = [fe[efx]] = [[fee] fx] = [ffx] = x$
and so $\delta$ is a surjection and, hence, an isomorphism. ∎

**Proposition 78.** *For* $e, f \in H$, $\delta : (H, \square_e) \cong (H, \square_f)$, *where* $\delta x = [xef]$.

Proof. We have $\delta (x \square_e y) = \delta [xye] = [[xye]ef] = [xy[eef]] = [xyf] = [[xee] yf] = [[xe[ffe]] yf] = [[[xef]fe] yf] =$
$= [[xef] [yef] f] = \delta x \square_f \delta y$ and so $\delta$ is a homomorphism. If $\delta x = \delta y$ then $[xef] = [yef]$ and so $x = [xee] =$
$= [xe [ffe]] = [[xef] fe] = [[yef] fe] = [ye [ffe]] = [yee] = y$. Therefore, $\delta$ is an injection. Since $\delta [xfe] = [[xfe] ef] =$
$= [xf [eef]] = [xff] = x$, $\delta$ is a surjection and, therefore, an isomorphism. ∎

Note 2 and Propositions 77 and 78 can be applied to prove the following Theorem. We omit the proof but we highlight that $\Omega\Psi(G) \cong G$, not equal to $G$ and $\Pi\Lambda(W) \cong W$, not equal to $W$.

**Theorem 79.** *Let* $G \in \mathbf{G}$, $W \in \mathbf{WQ}$ *and* $H \in \mathbf{H}$. *Then we have* $\Gamma\Phi(G) = G$, $\Phi\Gamma(W) = W$, $\Pi\Lambda(W) \cong W$, $\Lambda\Pi(H) = H$, $\Psi\Omega(H) = H$ *and* $\Omega\Psi(G) \cong G$.

**Theorem 80.** $[W]$ *is the natural ternary product of a Ward quasigroup* $(W, e)$ *if and only if*:

(80.1) $[xxe] = [yye] = e$ *for all* $x, y \in W$, *and for all* $x, y, z \in W$,
(80.2) $[xyz] = [[xye] ze] = [x [z[eye] e] e]$,
(80.3) $[xye] = [xze]$ *implies* $y = z$ *and*
(80.4) $[xye] = [wye]$ *implies* $x = w$.

Proof. ($\Rightarrow$) If $(W, e)$ is a Ward quasigroup with right unit $e$ then it follows from Proposition 75 that the natural ternary product $[W]$ satisfies (80.1) through (80.4)
($\Leftarrow$) If we define $xy = [xye]$ then, by Proposition 75, $W$ is a Ward quasigroup with right unit $e$. ∎

**Corollary 81.** *Let* $[W]$ *be the natural ternary product of a Ward quasigroup* $(W, e)$. *Then*:
(81.1) $e$ *is right unitary and*
(81.2) $e$ *is an outer lateral unit.*

Proof. Note that, by (80.1), $[eee] = e$. Then, by (80.2), $[[eee] xe] = [e[x[eee] e] e] = [e[xee] e]$ and so, by (80.3), $x = [xee]$ and $e$ is right unitary. Also, $[exe] = [e [e[exe] e] e]$ and, by (80.3), $x = [e[exe] e]$ and so $e$ is an outer lateral unit. ∎

**Corollary 82.** $[W]$ *is the natural ternary product of a Ward quasigroup* $(W, e)$ *if and only if*:

(82.1) *there exists a right unital element* $e \in W$,
(82.2) *for all* $x, y, z \in W$, $[xye] = [[xze] [yze] e]$ *and*
(82.3) *for all* $x, y \in W$ *there exists* $z \in W$ *such that* $x = [yze]$.

Proof. Necessity and sufficiency both follow from Proposition 75. ∎

**Definition 11.** $\mathbf{NWQ}$ is defined as the collection of all ternary products that are the natural ternary product of a Ward quasigroup.

**Theorem 83.** *There are mappings* $\Theta : \mathbf{G} \to \mathbf{NWQ}$ *and* $\Sigma : \mathbf{NWQ} \to \mathbf{G}$ *such that* $\Theta$ *and* $\Sigma$ *are inverse mappings.*

Proof. We define $\Theta(G,\bullet) = [xyz]_G$, where $[xyz]_G = z^{-1} \bullet y^{-1} \bullet x$. We define $\Sigma([W])$ as follows: $x \square_\Sigma y = [y\,[exe]\,e]$, where $[W]$ has the form as in Theorem 80. Then the fact that $\Sigma([W])$ is the dual of a group with identity element $e$ and $x^{-1} = [exe]$ follows from the proof of ($\Leftarrow$) in Theorem 80 and the comments regarding $(W,\circ_e)$ in Note 2. So $(\Sigma([W]), \square_\Sigma)$ is a group and $\Sigma$ is well defined. It is straightforward to check, using Theorem 80, that $\Theta(G) \in NWQ$. So $\Theta$ is well defined.

Consider $\Sigma \Theta(G)$. We have $x \square_\Sigma y = [y\,[exe]\,e] = e^{-1}\bullet[exe]^{-1}\bullet y = (x^{-1})^{-1}\bullet y = xy$ and so $\Sigma \Theta(G) = G$. Also,
$\Theta \Sigma([W]) = [xyz]_{\Sigma(S)} = z^{-1} \square_\Sigma y^{-1} \square_\Sigma x = [eze] \square_\Sigma [eye] \square_\Sigma x = ([eze] \square_\Sigma [eye]) \square_\Sigma x = [[eye]\,[e[eze]\,e]\,e] \square_\Sigma x =$
$= [[eye]\,ze] \square_\Sigma x$. Using Theorem 80, $[eyz] = [[eye]\,ze] = [e[z[eye]\,e]\,e]$ and so $[xyz]_{\Sigma(S)} = [[eye]\,ze] \square_\Sigma x =$
$= [eyz] \square_\Sigma x = [x[e[eyz]\,e]\,e] = [x\,[e\,[e[z[eye]\,e]\,e]\,e]\,e] = [x\,[z[eye]\,e]\,e] = [xyz]$. Therefore $\Theta \Sigma([W]) = [W]$. ∎

**Definition 12.** We define mappings $\alpha : NWQ \to H$ and $\beta : H \to NWQ$ as follows. For $[W] \in NWQ$, $\alpha[W] = [\,]_{[W]}$, where $[xyz]_{[W]} = [[xye][eze]e]$. Note that since $[W] \in NWQ$, by Theorem 80, such an element $e$ exists. We will prove below that $\alpha[W] \in H$. For $[H] \in H$ we define $\beta[H] = [\,]_h$, where $h \in H$ and $[xyz]_h = [[xyh]\,zh]$. We will show below that $[\,]_h \cong [\,]_k \in NWQ$.

**Theorem 84.** $\alpha[W] \in H$

Proof. Since $[W] \in NWQ$, Theorem 80 applies to $[W]$. So $[xxz]_W = [[xxe][eze]\,e] = [e[eze]\,e] = z$. Also, $[xyy]_{[W]} = [[xye][eye]\,e] = [xy\,[eye]] = [x\,[[eye][eye]e] = [xee] = x$. So every element is bi-unitary.

By the remark at the very start of this section, to prove that $\alpha[W] \in H$ we need only prove that $[[xyz]_{[W]}\,uv]_{[W]} = [xy\,[zuv]_{[W]}]_{[W]}$. Note that by corollary 82, we have
$\qquad$ (1): $[xye] = [[xze][yze]\,e] = [[xxe][yxe]\,e] = [e[yxe]\,e]$.
Now, using Theorem 80, $[[xyz]_{[W]}\,uv]_{[W]} = [\,[[[xye][eze]\,e]\,ue]\,[eve]\,e\,] = [\,[[xye][eze]\,u]\,[eve]\,e\,] =$
$= [\,[[xye]\,[u[eze]\,e]\,e]\,[eve]\,e\,] = [\,[[xye][uze]\,e]\,[eve]\,e\,] = [\,[xye]\,[\,[eve]\,[e[uze]\,e]\,e\,]\,]$. Using (1) twice,
$[[xyz]_{[W]}\,uv]_{[W]} = [\,[xye]\,[[eve]\,[e[uze]\,e]\,e]\,] = [\,[xye]\,[[eve]\,[zue]\,e]\,] = [\,[xye]\,[e\,[[zue][eve]\,e]\,e]\,] =$
$= [\,[xye]\,[e\,[zuv]_{[W]}\,e]\,] = [xy\,[zuv]_{[W]}]_{[W]}$. ∎

**Theorem 85.** $[\,]_h \cong [\,]_k \in NWQ$.

Proof. We use Theorem 80 and the fact that $[\,]$ is a heap to prove that $[\,]_h \in NWQ$. Recall that $[xyz]_h =$
$= [[xyh]\,zh]$. So $[xxh]_h = [[xxh]\,hh] = [xxh] = h$, which satisfies (80.1) when $h = e$. Now $[x\,[z[hyh]\,h]\,h] =$
$= [x\,[[zhy]\,hh]\,h] = [x\,[zhy]\,h] = [[xyh]\,zh] = [xyz]_h$, which satisfies (80.2). Then, $[xyh]_h = [xzh]_h$ implies $[[xyh]\,hh] = [[xzh]\,hh]$ implies $[xyh] = [xzh]$, which by Proposition 74 implies $y = z$. Similarly, $[xyh]_h = [wyh]_h$ implies $x = w$. So $[\,]_h$ satisfies (80.3) and (80.4). Hence, $[\,]_h \in NWQ$.

Now we prove that $\delta : [\,]_h \cong [\,]_k$, where $\delta x = [khx]$. By Proposition 74, $\delta$ is an injection. Since $\delta[hkx] =$
$= [kh\,[hkx]] = [[khh]\,kx] = [kkx] = x$, $\delta$ is a surjection and, hence, a bijection. So $[(\delta x)(\delta y)(\delta z)]_k =$
$= [[[khx][khy]\,k]\,[khz]\,k\,] = [[[khx]\,yh]\,[khz]\,k] = [\,[[khx]\,yh]\,zh\,] = [\,[kh\,[xyh]]\,zh\,] = [kh\,[[xyh]\,zh]] =$
$= \delta\,[[xyh]\,zh] = \delta\,[xyz]_h$. Thus, $\delta$ is a ternary product isomorphism; that is, $\delta : [\,]_h \cong [\,]_k$. ∎

**Theorem 86.** $\beta\,\alpha[W] \cong [W]$

Proof. $\beta \alpha [W] = \beta [\ ]_{[W]} = [\ ]_h = [[xyh]_{[W]} zh]_{[W]} = [[[xye][ehe] e] zh]_{[W]} = [\ [[[xye][ehe] e] ze] [ehe] e\ ]$.
But $[\ ]_h \cong [\ ]_e$. Then $[xyz]_e = [[xye] ze] = [xyz]$, so $\beta \alpha [W] = [\ ]_h \cong [\ ]_e = [\ ] = [W]$. ∎

**Theorem 87.** $\alpha \beta [H] = [H]$

Proof. $\alpha \beta [H] = \alpha [\ ]_h = [\ ]_{[W]}$ , where $[W] = [\ ]_h$. Now $[xyz]_{[W]} = [[xyh]_{[W]} [hzh]_{[W]} h]_{[W]}$ . But
$[xyh]_{[W]} = [[xyh]_h [hhh]_h h]_h = [xyh]_h = [[xyh] hh] = [xyh]$ and $[hzh]_{[W]} = [\ [hzh]_h [hhh]_h h\ ]_h =$
$= ([hzh]_h)_h = ([[hzh] hh])_h = [hzh]_h = [[hzh] hh] = [hzh]$. Therefore, $[xyz]_{[W]} = [[xyh]_{[W]} [hzh]_{[W]} h]_{[W]} =$
$= [[xyh][hzh] h]_h = [\ [[xyh][hzh] h] [hhh] h\ ] = [[xyh][hzh] h]$. But $[\ ]$ is a heap and so $[xyz]_{[W]} = [[xyh][hzh] h] =$
$= [[[xyh] hz] hh] = [[xyh] hz] = [xy [hhz]] = [xyz]$. So we have proved that $\alpha \beta [H] = [\ ]_{[W]} = [H]$. ∎

We can now extend Theorem 79 as follows:

**Theorem 88.** *Let* $G \in \mathbf{G}$ , $W \in \mathbf{WQ}$ , $[W] \in \mathbf{NWQ}$ *and* $H \in \mathbf{H}$ . *Then we have* $\Gamma\Phi(G) = G$, $\Phi\Gamma(W) = W$,
$\Pi\Lambda(W) \cong W$, $\Lambda\Pi(H) = H$, $\Psi\Omega(H) = H$, $\Omega\Psi(G) \cong G$, $\alpha\beta[H] = [H], \beta\alpha[W] \cong [W], \Sigma\Theta(G) = G$ *and*
$\Theta\Sigma([W]) = [W]$.

Although Theorem 88 allows us the move freely between the elements of $\mathbf{G}$ , $\mathbf{WQ}$ , $\mathbf{H}$ and $\mathbf{NWQ}$ , from one algebraic structure to another, it is not clear whether there are any advantages in doing so. For example, can a major, unproven conjecture in one algebraic class be shed light on, through a change of perspective to another class? Do results in one class hint at corresponding results in another class, results that perhaps have not even been considered?

## 10. Generalsied heaps that are standard or ternary operations on inverse semigroups

**Theorem 89.** *A ternary operation $[S]$ is the standard ternary operation of an inverse semigroup $S$*
*if and only if $[S]$ is a generalized heap and there exists an involutive unary mapping* $' : x \to x'$ *and an*
*idempotent unary mapping* $\wedge : x \to x^\wedge$ *on $S$ such that the following identities hold:*

(89.1) $[xx^\wedge y] = [xy'y^\wedge]$,
(89.2) $[x^\wedge x^\wedge x'] = x'$ and
(89.3) $[x'x'x^\wedge] = x^\wedge$.

Proof. ($\Rightarrow$) Suppose that a ternary operation $[S]$ is the standard ternary operation on an inverse semigroup $S$.
That is, $[xyz] = xy^{-1}z$. Then $[S]$ is a generalized heap (cf. observation 1 above). Define $x' = x^{-1}$ and $x^\wedge = x^{-1}x$.
Then $[xx^\wedge y] = x(x^{-1}x)^{-1}y = xy = xyy^{-1}y = x(y^{-1})^{-1}(y^{-1}y) = [xy'y^\wedge]$ , $x' = x^{-1} = (x^{-1}x)(x^{-1}x)^{-1}x^{-1} = [x^\wedge x^\wedge x']$
and $[x'x'x^\wedge] = x^{-1}(x^{-1})^{-1}(x^{-1}x) = x^{-1}x = x^\wedge$. Therefore, 89.1, 89.2 and 89.3 are valid.

($\Leftarrow$) Suppose that the ternary operation $[S]$ is a generalized heap with unary operations $'$ and $\wedge$, involutive and idempotent respectively, that satisfies 89.1 and 89.2. Define a binary operation $\sim$ on $S$ as follows:
(1): $x \sim y = [xx^\wedge y]$ for $x, y \in S$. We will denote $x \sim y$ by $xy$. Since $[S]$ is a generalized heap, $x(yz) = [xx^\wedge[yy^\wedge z]] =$
$= [[xx^\wedge]y^\wedge z] = [[xy'y^\wedge]y^\wedge z] = [x[y^\wedge y^\wedge y']z] = [xy'z]$. That is,

$$(2): x(yz) = [xx^\wedge[yy^\wedge z]] = [xy'z].$$

Then, using 89.1, $(xy)z = [xx^\wedge y]z = [xy'y^\wedge]z = [[xy'y^\wedge][xy'y^\wedge]^\wedge z] = [[xy'y^\wedge]z'z^\wedge] = [xy'[y^\wedge z'z^\wedge]] =$
$= [xy'(y^\wedge z)] = [xy'[y^\wedge y^{\wedge\wedge} z]] = [xy'[y^\wedge y^\wedge z]]$ . So

$$(3): (xy)z = [xy'[y^\wedge y^\wedge z]].$$

But then (2) and (3) imply $(xy)z = [xy\,'\,[y^\wedge y^\wedge z]] = x(y[y^\wedge y^\wedge z]) = x\,[yy^\wedge[y^\wedge y^\wedge z]] = x[y[y^\wedge y^\wedge y^\wedge]z] = = x[yy^\wedge z] = [xx^\wedge[yy^\wedge z]] = x(yz)$. So the product $xy = [xx^\wedge y]$ gives a semigroup structure to $S$. Now $S$ is a regular semigroup since, by 2, $xx\,'\,x = [xx\,'\,'\,x] = [xxx] = x$. We need to show that the idempotents of $S$ commute, for then, by Theorem 1.17 [1], $\{S,\sim\}$ is an inverse semigroup.

First we show that for any $y,z \in S$, $y\,'\,y$ and $z\,'\,z$ commute. Since $[S]$ is a generalized heap, for any $x,y,z \in S$, we have $[[xyy]zz] = [[xzz]yy]$. By (2), $[xyz] = xy\,'\,z$ and so

(4): $xy\,'\,yz\,'\,z = xz\,'\,zy\,'\,y$  and, setting $x = y\,'\,y$ gives $y\,'\,yy\,'\,yz\,'\,z = y\,'\,yz\,'\,zy\,'\,y$ and, since $y\,'\,yy\,'\,y = y\,'\,y$, we have that (5): $y\,'\,yz\,'\,z = y\,'\,yz\,'\,zy\,'\,y$. Similarly, the identity $[yy[zzx]] = [zz[yyx]]$ implies that (6): $y\,'\,yz\,'\,z = z\,'\,zy\,'\,yz\,'\,z$.

Since $y\,'\,y$ is idempotent, (6) implies that $y\,'\,yz\,'\,z$ is idempotent. Hence, (5) and (6) imply

(7): $y\,'\,yz\,'\,z = y\,'\,yz\,'\,zy\,'\,y = z\,'\,zy\,'\,yz\,'\,zy\,'\,y = z\,'\,zy\,'\,y$. In particular, (8) $xx\,'\,x\,'\,x = x\,'\,xxx\,'$. We now show that every idempotent element $x$ is of the form $x = xx\,'$. We have $x = xx$ and so by (8), $xx\,'\,x\,'\,x = x\,'\,xxx\,' = x\,'\,xx\,' = x\,'$. Hence, $x\,'\,x = (xx\,'\,x\,'\,x)x = xx\,'\,x\,'\,x = x\,'$ and then $xx\,' = xx\,'\,x = x$. Then, since (7) implies that $xx\,'$ and $yy\,'$ commute, the idempotents of $S$ commute and, as previously noted, $S$ is an inverse semigroup with $x\,'$ equal to the inverse of $x$. Therefore, since (2) implies $[xyz] = xy\,'\,z$, $[S]$ is the standard ternary operation of the inverse semigroup $\{S,\sim\}$. ∎

**Corollary 90.** *Standard ternary operations of inverse semigroups form an algebraic variety of type* $(3,1,1)$.

Note that in the proof of Theorem 89 ($\Leftarrow$) the condition 89.3 is not used. However it is needed in order to prove Theorem 91 below.

**Definition.** We define *GHI* to be the collection of generalized heaps that appear as the standard ternary operation of an inverse semigroup and *I* to be the collection of inverse semigroups. We define mappings $g:GHI \to I$ and $i:I \to GHI$ as follows: $g([S],\,',\,^\wedge) = \{S,\sim\}$, where $\sim$ is defined as in the proof of Theorem 89 ($\Rightarrow$), and $iS = ([S]\text{-}1,\,^{-1},\,^\wedge)$, where $x^\wedge = x^{-1}x$.

**Theorem 91.** *The mappings g and i are inverse mappings.*

The proof of Theorem 91 is straightforward and is omitted. However, we note that condition 89.3 implies, by definition, that $x^\wedge = x\,'\,x = x^{-1}x$. Without condition 89.3 there is no apparent way to prove that $x^\wedge = x^{-1}x$, which is necessary in order to prove that $ig([S],\,',\,^\wedge) = i\{S,\sim\} = ([S],\,',\,^\wedge)$.

**Theorem 92.** *A ternary operation $[S]$ is the natural ternary operation of an inverse semigroup $S$ if and only if $[S]$ is associative and there exists an involutive unary mapping $\,': x \to x\,'$ and an idempotent unary mapping $\,^\wedge: x \to x^\wedge$ on $S$ such that the following identities hold:*

(92.1) $x = [xx\,'\,x] = [xx^\wedge x^\wedge]$,
(92.2) $[xx^\wedge y] = [xyy^\wedge]$,
(92.3) $x\,' = [x^\wedge x^\wedge x\,']$,
(92.4) $x^\wedge = [x\,'\,xx^\wedge] = [x^\wedge x\,'\,x]$,
(92.5) $[[xx^\wedge x\,'\,][y\,'\,yy^\wedge][y\,'\,yy^\wedge]^\wedge] = [[y\,'\,yy^\wedge][xx^\wedge x\,'\,][xx^\wedge x\,'\,]^\wedge]$ *and*
(92.6) $[xyz] = [[xx^\wedge y]zz^\wedge]$.

Proof ($\Rightarrow$) Suppose that $[S]$ is the natural ternary operation of an inverse semigroup $S$. Then since the semigroup $S$ is associative, it is clear that $[S]$ is associative. If we define $\,' = ^{-1}$ and $x^\wedge = x^{-1}x$ then it is straightforward to prove that $\,'$ and $^\wedge$ are involutive and idempotent unary operations respectively and that the identities 92.1, 92.2, 92.3, 92.4 and 92.6 are satisfied. Then the identity 92.5 implies that $xx^{-1}y^{-1}y = y^{-1}yxx^{-1}$, which is valid because the idempotents of an inverse semigroup commute.

($\Leftarrow$) Suppose that [S] is an associative ternary operation that satisfies the identities 92.1, 92.2, 92.3, 92.4, 92.5 and 92.6. Define a binary operation ~ on S as follows: (1): $x$~$y$ = $[xx^\wedge y]$ for $x,y \in S$. We will denote $x$~$y$ by $xy$. Then (1), 92.2 and the associativity of [S] imply that $(xy)z = [[xx^\wedge y]zz^\wedge] = [xx^\wedge [yzz^\wedge]] = x(yz)$.

Also, (1), 92.1, 92.4 and associativity imply that $(xx')x = [[xx^\wedge x']xx^\wedge] = [xx^\wedge [x'xx^\wedge]] = [xx^\wedge x^\wedge] = x$. So {$S$,~} is a regular semigroup. As noted in the proof of Theorem 89 ($\Leftarrow$), a regular semigroup with commuting idempotents is an inverse semigroup. So we only need to show that the idempotents of {$S$,~} commute, since by 92.6, $xyz = [xyz]$.

Firstly, {$S$,~} is a semigroup, $(xx')(xx') = \{(xx')x\}x' = xx'$. So $xx'$ is idempotent. We proceed to show that every idempotent $x$ satisfies $x = xx'$. Suppose that $x = xx$. By (1) and 92.5 $xx'y'y = y'yxx'$ and so $xx'x'x = x'x xx' = x'xx' = x'$. Then $x'x = (xx'x'x)x = xx'x'x = x'$ and $xx' = x(x'x) = x$. So any two idempotents $x$ and $z$ can be expressed as $x = xx'$ and $z = zz'$. Since we know that $xx'y'y = y'yxx'$, setting $y = z'$ gives $xx'zz' = zz'xx'$, or, $xz = zx$. ∎

We are now in a position to expand the table that appears on page 1. Note that the table below distinguishes between the natural ternary product [S] and the ternary product [S]∗ determined by the unary operation ∗.

| *Ternary Product Type* {S}= [S] | *Ternary Product Type* {S} = [S]∗ | | *Groupoid Type* |
|---|---|---|---|
| | l-bi-unital semiheap (∗ = ') | | Involuted monoid with identity $l$ and involution ' |
| | Idempotent, l-bi-unital semiheap (∗ = ') | | Monoid with identity $l$, involution ' and $xx'x = x$ |
| | l-bi-unital generalised heap (∗ = -1) | | Inverse monoid with identity $l$ |
| | l-bi-unital heap (∗ = -1) | | Group with identity $l$ |
| Laterally commutative, l-bi-unital semiheap | | | Right modular with left identity $l$ |
| | Left and Laterally commutative, l-unital, $\{xy^\wedge l\}=\{x^\wedge yl\}$ ^ is idempotent (∗ = ') | | Right modular with left identity $l$ |
| l-bi-unital, left l-consistent, l central commutant, ^ -congruent, $\{lxl\} = x^\wedge$, ^ is idempotent | | | ∗-unary with left identity |

| *Ternary Product Type* | *Ternary Product Type* | | *Groupoid Type* |
|---|---|---|---|
| $\{S\} = [S]$ | $\{S\} = [S]*$ | | |
| | $l$-unital, left commutative, left $l$-consistent, $\{xly\}=\{ylx\}$, $\{xyz\}^\wedge = \{x^\wedge y^\wedge z^\wedge\}$, $^\wedge$ is idempotent, $l = l^\wedge$ | | $*$-unary with left identity $l$ |
| Right commutative, $l$ right unital, $\{x\{yzq\}w\} = \{zw\{yxq\}\}$, $\{\{xyz\}qw\} = \{\{xyq\}wz\}$ | | | The dual of a right modular groupoid with left identity $l$ |
| $e$ right unital, $\{xye\} = \{\{xze\}\{yze\}e\}$, For all $x$ and $y$ there exists $z$ such that $x = \{yze\}$ | | | A Ward quasigroup; that is, $xx = yy = e$ ; $x(yy) = x$ ; $xy = (xz)(yz)$ |
| | Generalised heap, $\{xx^\wedge y\} = \{xy\,'y^\wedge\}$ $\{x^\wedge x^\wedge x\,'\} = x\,'$ $\{x\,'x\,'x^\wedge\} = x^\wedge$ $^\wedge$ is idempotent $\,'$ is involutive $(*= -1)$ | | An inverse semigroup |
| Associative $x = \{xx\,'x\} = \{xx^\wedge x^\wedge\}$ $\{xx^\wedge y\} = \{xyy^\wedge\}$ $x\,' = \{x^\wedge x^\wedge x\,'\}$ $x^\wedge = \{x\,'xx^\wedge\} = \{x^\wedge x\,'x\}$ $\{\{xx^\wedge x\,'\}\{y\,'yy^\wedge\}\{y\,'yy^\wedge\}^\wedge\}=$ $\{\{y\,'yy^\wedge\}\{xx^\wedge x\,'\}\{xx^\wedge x\,'\}^\wedge\}$ $\{xyz\} = \{\{xx^\wedge y\}zz^\wedge\}$ $^\wedge$ is idempotent $\,'$ is involutive | | | An inverse semigroup |

**Example 3.** The following is an example of a generalised heap that is *not* the standard ternary operation of an inverse semigroup. The author is grateful to Professor Boris M. Schein for this example.

We have $S = \{0,a,b\}$ with $[aaa] = a = [abb]$, $[bbb] = b = [baa]$ and all other ternary products equal to 0. Routine calculations show that $\{S,[\ ]\}$ is a generalised heap. Assume that $[\ ]$ is the *standard* ternary operation of an inverse semigroup $S$. Then $[xyz]\,' = [z\,'y\,'x\,']$, where $\,' = -1$.

Now, since $\,'$ is an involutive mapping on $S$, either $\,'=1_S$, $0\,' = a$ and $b\,' = b$, $0\,' = b$ and $a\,' = a$ or $a\,' = b$ and $0\,' = 0$. In the first case, $a = [abb] = [abb]\,' = [b\,'b\,'a\,'] = [bba] = 0$, a contradiction. In the second case, $b = [baa] = [baa]\,' = [a\,'a\,'b\,'] = [00b] = 0$, a contradiction. In the third case, $a = [abb] = [abb]\,' = [b\,'b\,'a\,'] = [00a] = 0$, a contradiction. In the fourth and last case, $a = b\,' = [baa]\,' = [a\,'a\,'b\,'] = [bba] = 0$, a contradiction.

**Example 4.** If *S* = {0,a,b} and the ternary operation < > satisfies the identity < *xyz*> = *x* then the associative ternary operation {*S*,< >} is ***not*** the natural ternary operation on an inverse semigroup *S*. This is because *x*' = [*xyz*]' = [*z*' *y*' *x*'] = *z*', where ' = -1, a contradiction.

10 Albert Mansions, Crouch Hill, London, N8 9RE, United Kingdom
bobmonzo@talktalk.net
00442083406167